\documentclass{article}
\usepackage{amsmath,amssymb,enumerate,bbm,calc,capt-of,ifthen}
\usepackage{graphicx}% Include figure files
\usepackage{subfigure}

\usepackage{bm}% bold math
\usepackage{algorithm2e}
\dontprintsemicolon
\SetKw{klet}{let}
\SetCommentSty{textit}
\SetKwFunction{ifft}{InverseFFT}
\linesnumbered
\restylealgo{algoruled}

\usepackage{epsfig}

\newtheorem{cntr}{do not use}

\newtheorem{definition}[cntr]{Definition}
\newtheorem{assumption}{Assumption}

\newtheorem{varremark}[cntr]{Remark}
\newenvironment{remark}{\begin{varremark}\em}{\em\end{varremark}}

\newtheorem{varnote}{Note}

\newtheorem{proposition}[cntr]{Proposition}

\newenvironment{proof}{
 \noindent\textbf{Proof.}\ }{\hspace*{\fill}
 \begin{math}\Box\end{math}\medskip}

\newenvironment{proofof}[1]{
  \noindent\textbf{Proof of #1.}\ }{\hspace*{\fill}
  \begin{math}\Box\end{math}\medskip}

\newenvironment{proof*}[1]{
  \noindent\textbf{#1\ }}{\hspace*{\fill}
  \begin{math}\Box\end{math}\medskip}

\newtheorem{lemma}[cntr]{Lemma}
\newtheorem{corollary}[cntr]{Corollary}
\newtheorem{theorem}[cntr]{Theorem}

\numberwithin{cntr}{section}
\numberwithin{equation}{section}

\newcommand{\comment}[1]{}
\newcommand{\abs}[1]{\left| #1 \right|}%Absolute value
\newcommand{\absSmall}[1]{| #1 |}%Absolute value, small
\newcommand{\RR}[0]{{\mathbb{R}}}
\newcommand{\Sone}[0]{{\mathbb{S}^{1}}}

\newcommand{\vk}[0]{{\vec{k}}}

\newcommand{\va}[0]{{\vec{a}}}
\newcommand{\vb}[0]{{\vec{b}}}

\newcommand{\vp}[0]{{\vec{p}}}

\newcommand{\Oto}[1]{{0 \ldots #1-1}}

\newcommand{\curvemax}{\bar{\kappa}}
\newcommand{\curvemin}{\underline{\kappa}}
\newcommand{\kmax}{{k_{\textrm{max}}}}
\newcommand{\ktex}{{k_{\textrm{tex}}}}
\newcommand{\curvemaxi}{{\curvemax^{-1}}}

\newcommand{\arcs}[2]{{\mathfrak{A}_{#1,#2}}}
\newcommand{\arcsMinusOne}[3]{{\arcs{#1}{#2}^{#3}}}

\newcommand{\kp}{{{\vk}^{\perp}}}
\newcommand{\kr}{{k_{r}}}
\newcommand{\kt}{{k_{\theta}}}
\newcommand{\ktv}{{\vec{k}_{\theta}}}

\newcommand{\curvesep}{{\delta}}
\newcommand{\pointNoise}{{\zeta}}
\newcommand{\tanNoise}{{\xi}}

\newcommand{\step}[1]{{1_{\gamma_{#1}}}}
\newcommand{\fstep}[1]{{\widehat{1_{\gamma_{#1}}}}}

\newcommand{\screen}{{[0,1]^2}}

\newcommand{\ZZf}[1]{{\mathbb{Z}_{2^{#1}}}}

\newcommand{\Img}{{\bm{\rho}}}
\newcommand{\ImgTex}{{\Img_{\textrm{tex}}}}
\newcommand{\Imgk}{\widehat{\Img}}
\newcommand{\ImgCoeff}[1]{{\Img_{#1}}}
\newcommand{\contrastMin}{{\underline{\Img}}}
\newcommand{\contrastMax}{{\overline{\Img}}}

\newcommand{\RadFilt}{{\mathcal{W}}}
\newcommand{\RadFiltCentered}{{\mathcal{W}^{p}}}
\newcommand{\RadFiltCenteredX}{{\check{\RadFiltCentered}}}
\newcommand{\RadFiltX}{{\check{\RadFilt}}}
\newcommand{\RadFiltXSup}{{\overline{\RadFilt}}}
\newcommand{\AngFilt}{{\mathcal{V}}}

\newcommand{\ConstRadFiltDecay}{{\mathfrak{C}_{\RadFilt}}}
\newcommand{\ConstTangentNormal}{{\mathfrak{C}(\RadFilt, \AngFilt,\alpha)}}
\newcommand{\SurfelResolution}{{\mathfrak{D}(\RadFilt, \AngFilt, \alpha)}}
\newcommand{\SurfelThreshold}{{\mathfrak{T}(\RadFilt, \AngFilt, \alpha)}}
\newcommand{\ConstGeo}{{\mathfrak{C}_{\textrm{geo}}}}

\newcommand{\norm}[2]{\left\| #1 \right \|_{#2}}%Compute the norm

\newcommand{\fourier}{{\mathcal{F}}}
\newcommand{\ifourier}{{\mathcal{F}^{-1}}}
\newcommand{\dfilter}[2]{{D_{#1,#2}}}

\begin{document}

\title{Spectral edge detection in two dimensions using wavefronts}

\author{L. Greengard and C. Stucchio}

\maketitle

\begin{abstract}
A recurring task in image processing, approximation theory, and the
numerical solution of partial differential equations is to
reconstruct a piecewise-smooth real-valued function
$f(x)$, where $x \in \RR^N$,
from its truncated Fourier transform
(its truncated spectrum). An essential step is {\em edge detection}
for which a variety of one-dimensional schemes
have been developed over the last few decades. Most higher-dimensional
edge detection algorithms consist of applying
one-dimensional detectors in each component direction in order to recover the
locations in $\RR^N$ where $f(x)$ is singular (the singular support).

In this paper, we present a multidimensional algorithm which identifies
the {\em wavefront} of a function from spectral data.
The wavefront of $f$ is the set of points
$(x,\vk) \in \RR^{N} \times (S^{N-1} / \{\pm 1\})$ which encode
both the location of the singular points of a function and the orientation
of the singularities. (Here $S^{N-1}$ denotes the unit sphere in $N$ dimensions.) More precisely, $\vk$ is the direction of the
normal line to the curve or surface of discontinuity at $x$.
Note that the singular support is simply the projection of the
wavefront onto its $x$-component.
In one dimension, the wavefront is a subset of
$\RR^{1} \times (S^{0} / \{\pm 1\}) = \RR$, and it coincides with the singular
support. In higher dimensions, geometry comes into play and they are distinct.
We discuss the advantages of wavefront reconstruction and indicate how
it can be used for segmentation in magnetic resonance imaging
(MRI).
\end{abstract}

\section{Introduction}

Consider an image, i.e. a function $\Img: \RR^{2} \rightarrow \RR^{+}$. If the image is smooth (e.g. $C^{m}(\RR^{2})$), then the Fourier transform of $\Img(x)$,
denoted $\Imgk(\vk)$, will decay rapidly (and hence be localized near $\vk=0$).
Discontinuities in the image cause $\Imgk(\vk)$ to decay more
slowly as $|\vk| \rightarrow \infty$.
Thus, information about the discontinuities can be said to be encoded in
the high frequency components of $\Imgk(\vk)$.
The goal of spectral edge detection is to recover the location of
the discontinuities  from limited (and often noisy) information about
$\Imgk(\vk)$.

As an example, consider the function $\Img(x) = 1_{B}(x)$ which is equal to $1$ inside $B = \{ x : \abs{x} < 1 \}$ and $0$ elsewhere.
The set of discontinuities of this function (the singular support) is given
by $\{ x : \abs{x} = 1 \}$.
One natural approach to computing the curve on which the discontinuities lie
is to first find a set of point which lie in the
singular support, followed by an algorithm aimed at connecting these
points sets into a finite number of curves.
In our example, the output would be the unit circle.
A relatively recent and important class of methods for locating
the singular support is based on concentration kernels
\cite{tadmor:edgesReview, gelb:gegenbauerReconstructionError,gelbcates:mriSegmentationSpectralData,gelbtadmor:edges,588632},
a high pass filtering approach, which we describe briefly
below.

For the function $\Img(x) = 1_{B}(x)$, the normal at each
point of discontinuity $x$ is simply the normal to the unit circle
at that point.  The \emph{wavefront} of this function is,
%%therefore, $\{ (x, k) : \abs{x} = 1 \textrm{~and~} \vk \parallel x \}$.
therefore, $\{ (x, k) \}$, with $\abs{x} = 1$ and $\vk \parallel x$.
In this paper, we study the problem of extracting the wavefront from
continuous spectral data available in a finite frequency range
$|\vk| < k_{max}$ in two dimensions.
This extra information is useful practically as well as theoretically.
First, it is easier to reconstruct curves of discontinuity from
points in the wavefront than points in the singular support, both in
closing ``gaps'' and in associating points on close-to-touching curves
to the correct one \cite{amenta98crust,chengnoise,dey99curve,dey01reconstructing,meleslie:curveppaer}.
Second, the directional information is useful for noise
removal. If spurious points are included in the
wavefront set, the normal (or tangent) data allows us to filter it out;
it is unlikely that a random point and a random tangent will be consistent
with the points and tangents that come from the actual curves of
discontinuity, as we have shown previously in
\cite{meleslie:curveppaer}.

Our approach to edge detection is based on applying concentration
kernels (high pass filters) to angular slices of the Fourier data.
Rather than recovering the points on the edges
(as in \cite{tadmor:edgesReview,gelbtadmor:edges,588632}),
we also determine the direction of the normal.
In section \ref{sec:mathformulation}, we present a precise
statement of the problem to be solved and
an overview of the full edge detection procedure (section \ref{sec:informal}).
In section \ref{sec:largeKAsymptotics}, we  present a detailed analysis
of the asymptotic behavior of the Fourier transform of the characteristic
function of a smooth region. Section \ref{sec:directionalFilters} is
devoted to a discussion of the {\em directional filters} used to extract
wavefront data and Section \ref{sec:connectingsurfels} describes the full algorithm.
Many of the proofs are technical and we
have relegated most proofs to the appendices.
We discuss the application of our method to
magnetic resonance imaging (MRI), where raw data is acquired
in the Fourier domain, extending the algorithm
of \cite{gelbcates:mriSegmentationSpectralData}).
Finally, we note that our algorithm is closely related
to the recent paper \cite{guolabatelim:edges}, which describes
a wavefront extraction procedure based on the {\em shearlet} transform,
extending the use of curvelets in \cite{MR2165380,MR2012649}.

\section{Mathematical Formulation}
\label{sec:mathformulation}

Let $\Img(x)$ be a piecewise smooth function form $\RR^2 \rightarrow \RR$,
and let $\chi(x)$ be a $C^{\infty}$ compactly supported radial
function. The \emph{singular support} of $\Img(x)$ consists of the
points $x_{0} \in \RR^{2}$ for which $\chi(\lambda(x - x_{0})) f(x)$ has
slowly decaying Fourier transform for every $\lambda > 0$, i.e.:
\begin{equation}
  \forall \lambda > 0, \sup_{\abs{\vk} \geq \kr} \abs{\int e^{i \vk \cdot x} \Img(x) \chi((x-x_{0})\lambda) dx } = O(\kr^{-3/2})
\end{equation}

In order to find the singular support,
concentration kernel methods
\cite{tadmor:edgesReview, gelb:gegenbauerReconstructionError,gelbcates:mriSegmentationSpectralData,gelbtadmor:edges,588632}
multiply the Fourier data $\Imgk(\vk)$ by a function which gives heavier weight to high frequencies than to low frequencies (a {\em high-pass} filter).
Since high frequencies encode the location of singularities but are unaffected by smooth parts of the image, this method isolates discontinuities from the rest of the image.
In short, concentration kernel methods find the \emph{location} of
singularities by flagging local maxima in the inverse Fourier transform  of
the high-pass filtered Fourier data.

The \emph{wavefront} of a function consists of the points $(x_{0}, \vk_{0}) \in \RR^{2} \times \Sone$ for which the Fourier transform of $\chi(\lambda(x - x_{0})) \Img(x)$ decays slowly in the direction $\vk_{0} = (\kr,\kt)$:
\begin{equation}
  \forall \lambda > 0, \sup_{r \geq \kr} \abs{\int e^{i r k_{0} \cdot x} \Img(x) \chi((x-x_{0})\lambda) dx } = O(\kr^{-3/2})
\end{equation}
As indicated in the introduction, while the singular support of $\Img(x)$
only contains the location of singularities, the wavefront also contains the
direction of the singularities.

\begin{remark}
In the language of computational geometry, the singular support
is a set of points, while the wavefront is a set of surfels
(pairs of the form $(x, k)$ with $x$ representing a
position and $k$ a direction).
\end{remark}

\subsection{Definition of the image class}

To simplify the theory, we consider a special class of images.
In particular, we consider two-dimensional images supported on $\screen$
and vanishing near the boundaries, which consist of a set of
piecewise constant functions on which is superimposed a globally smooth
function:
\label{eq:assumptions}
\begin{equation}
  \label{eq:imageClass}
  \Img(x) = \left[\sum_{j=0}^{M-1} \ImgCoeff{j}\step{j}(x) \right] + \ImgTex(x)
\end{equation}
where $\gamma_{j}(t)$ are simple closed curves, and $\step{j}(x)=1$ for $x$ in the interior of $\gamma_{j}$ and $0$ elsewhere. The ``texture'' term $\ImgTex(x)$ is band limited, i.e. $\widehat{\ImgTex}(\vk) = 0$ for $\abs{\vk} \geq \ktex$.

\begin{definition}
Let $\gamma_{j} (t)$ be a simple closed curve. The curvature at each point is denoted by
$\kappa_{j}(t)$, the normal to $\gamma_j(t)$ is denoted by $N_{j}(t)$, and the tangent is denoted
by $T_j(t)$.
\end{definition}

\begin{figure}
\setlength{\unitlength}{0.240900pt}
\ifx\plotpoint\undefined\newsavebox{\plotpoint}\fi
\sbox{\plotpoint}{\rule[-0.200pt]{0.400pt}{0.400pt}}%
\includegraphics[scale=0.5]{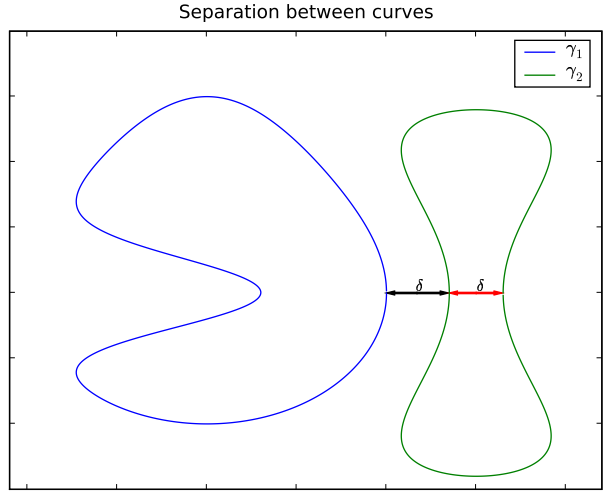}
\caption{An illustration of Assumption \ref{ass:separatedCurves}. The black arrow illustrates the condition \eqref{eq:separationAssumption}, while the red arrow illustrates the condition \eqref{eq:separationAssumptionSameCurve}.}
\label{fig:separationBetweenCurves}
\end{figure}

\begin{assumption}
  \label{ass:boundedCurvature}
  We assume the curves have bounded curvature, i.e.:
  \begin{equation}
    \label{eq:curvatureCondition}
    \forall i = \Oto{M}, ~ \kappa_{i}(s) = \frac{
      \abs{\gamma_{i,x}'(t) \gamma_{i,y}''(t) - \gamma_{i,y}'(t) \gamma_{i,x}''(t)}
    } {
      (\gamma_{i,x}'(t)^{2}+\gamma_{i,y}'(t)^{2})^{3/2}
    } \leq \curvemax
  \end{equation}
\end{assumption}
\begin{assumption}
  \label{ass:separatedCurves}
  We also assume that the curves are separated from each other (see Fig. \ref{fig:separationBetweenCurves}):
  \begin{subequations}
    \begin{equation}
      \label{eq:separationAssumption}
      \sup_{t,t'} \abs{ \gamma_{i}(t) - \gamma_{j}(t')} \geq \curvesep \textrm{~for~} i \neq j
    \end{equation}
    and that different areas of the same curve are separated from each
    other, i.e. if the curves are parameterized to move at unit speed,
    then:
    \begin{equation}
      \label{eq:separationAssumptionSameCurve}
      \sup_{\abs{t-t'} > \curvemaxi\pi/2 } \abs{ \gamma_{i}(t) - \gamma_{i}(t')} \geq \curvesep
    \end{equation}
  \end{subequations}
\end{assumption}

In order to extract the wavefront, we assume each edge is associated with a sufficiently large discontinuity. Since our wavefront detector decays somewhat slowly away from the wavefront,
we also assume that the discontinuity is bounded from above, so as not to pollute nearby
edges of lower contrast. We formalize this as:

\begin{assumption}
  \label{ass:contrastAssumption}
  We assume the contrast of the discontinuities in the image is bounded above and below:
  \begin{equation}
    \label{eq:contrastAssumption}
    0 < \contrastMin \leq \ImgCoeff{j} \leq \contrastMax
  \end{equation}
\end{assumption}

The data we are given and wish to segment are noisy samples of $\Imgk(\vk)$
obtained on a $2^{m} \times 2^{m}$ grid with spacing $2 \pi$ in $\vk$-space
centered at the origin:
$\vk \in 2\pi \ZZf{m}$ with $\ZZf{m}= \{ -2^{m-1},\ldots, 2^{m-1}-1 \}^{2}$.
Our segmentation goal is to recover the curves $\gamma_{j}(t)$.

With the preceding sampling, we may define the maximum frequency available in the image as
$\kmax = 2\pi2^{m-1}$. We assume that
$\ktex < \kmax$. More precisely, we assume that  $\kmax - \ktex \geq 12 \cdot 2\pi$, providing at least 12 lattice points in the sampling beyond $\ktex$. In most of our experiments, we take $\ktex=2\pi \cdot 16$, and $\kmax = 2\pi \cdot 32$.

Finally, we make the technical assumption that the curves have non-vanishing curvature:

\begin{assumption}
  \label{ass:curvatureBoundedBelow}
  We assume that the curvature of the curves is bounded below.
  \begin{equation}
    \label{eq:curvatureBoundedBelow}
    \forall i = \Oto{M}, ~ \kappa_{i}(s) \geq \curvemin > 0
  \end{equation}
\end{assumption}

\begin{remark}
  The assumption \eqref{eq:curvatureBoundedBelow} implies that the region bounded by $\gamma_{i}(t)$ is convex. This geometric fact is not used by our algorithm in any way.
\end{remark}

\begin{remark}
Assumption \ref{ass:curvatureBoundedBelow} is introduced only because it is
required below by our proof technique for the correctness of the directional filters. It is not strictly necessary, and it would be straightforward to extend our analysis to cases where the curvature vanishes. This, however, would require more complicated conditions on higher derivatives
of $\gamma_{j}(t)$ in places where the curvature vanishes, and correspondingly more complicated proofs. See Figure \ref{fig:squareExample} which illustrates that our algorithm works even when Assumption \ref{ass:curvatureBoundedBelow} is violated.
\end{remark}

\subsection{Wavefront Extraction Methodology}

The algorithm we present in this work takes an image, given as spectral data $\Imgk(\vk)$, and extracts a set of surfels sampled from the wavefront. We do this in two steps.

First, we construct a set of directional filters:
\begin{equation}
  \label{eq:filterDef}
  [\dfilter{\theta}{\alpha} \hat{\Img}](x) = \ifourier [\RadFilt(\kr) \AngFilt(\kt - \theta) \Imgk(\vk)]
\end{equation}
Here, $\RadFilt(\kr)$ is a high pass filter in the radial direction (the \emph{radial filter}), $\AngFilt$ is a smooth function, compactly supported on $[-\alpha,\alpha]$ (the \emph{angular filter}), and $\fourier$ is the Fourier transform. Note that, given a direction $\theta$, the angular filter is supported on the angular window $[\theta-\alpha,\theta+\alpha]$. The angular and radial filters will be related by the parabolic scaling:
\begin{equation*}
  \textrm{width}^{2} = \textrm{length}
\end{equation*}
(Equivalently, they will be chosen to satisfy $\textrm{width} = \textrm{length}^{2}$ in the $x$-domain.) In particular, this implies that the filter angle $\alpha \sim (\textrm{pass-band})^{-1/2}$, where
$\textrm{pass-band}$ is the center of the pass-band of the radial filter (to be determined
in section \ref{sec:directionalFilters}).

When applied to the function $\Imgk(\vk)$, the directional filters will return a function which is small except near locations where $N_{j}(t)$ points in the direction $\theta$. This allows us to pinpoint the locations of singularities with direction $\theta$. The result is, in some sense, a directional version of the jump function of \cite{tadmor:edgesReview,gelbcates:mriSegmentationSpectralData,gelbtadmor:edges,588632}. Spikes (local maxima) obtained from these directional filters correspond to {\em surfels} in the wavefront of the image.

\begin{remark}
For the algorithm described here to work, it is required that accurate values of the \emph{continuous} spectral data be available. The discrete Fourier transform (DFT) of $\Img(x)$ can not be used as a substitute for the continuous data, since aliasing induced by the DFT will destroy the asymptotic expansion of $\Imgk(\vk)$.
\end{remark}

\subsection{A note on our phantom}
\label{sec:PhantomExplanation}

\begin{figure}
  \includegraphics[scale=0.75, bb = 18pt 180pt 594pt 612pt]{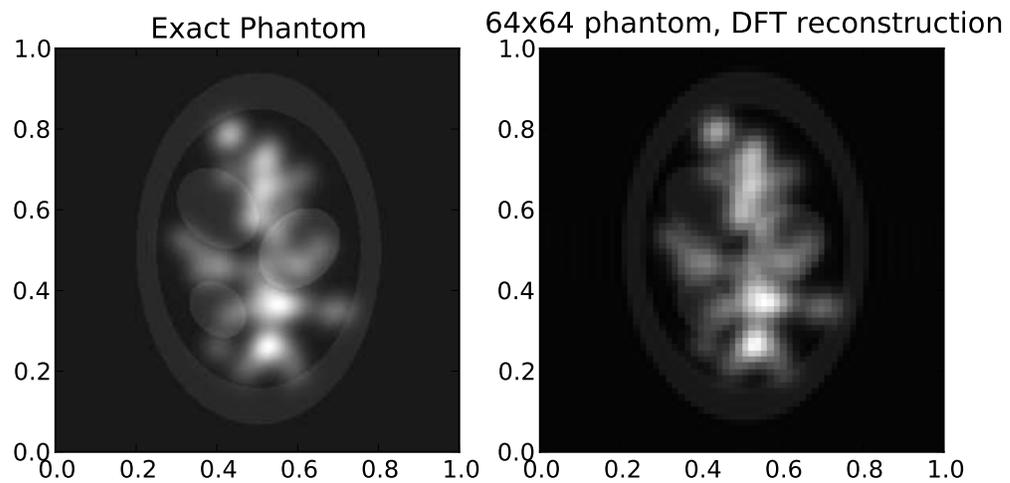}
  \caption{A plot of the phantom used in this work. The image in the second panel is obtained by applying the DFT to the (exactly known) continuous Fourier Transform of the phantom on a $64 \times 64$ grid of samples in $k$-space.}
  \label{fig:phantom}
\end{figure}

In the medical imaging literature, the Shepp-Logan phantom (which is
piecewise constant) is a traditional choice for analysis and validation
purposes. We have added a smoothly varying component
(see (\ref{eq:imageClass})), with
$\ImgTex(x)$ defined as a sum of Gaussians whose bandwidth (to six digits
of accuracy) is equal to half the bandwidth of the measurement data.

\subsection{Informal Description of the Algorithm}
\label{sec:informal}

The full algorithm proceeds in three steps.

\begin{enumerate}
\item Construct directional filters, based on
$\kmax$ (the maximum frequency content of the data) and
$\ktex$ (the frequency content of the ``smooth'' part of the image).
\item Multiply the data $\Imgk$ by each directional filter
indexed by $\theta = j\pi/A,\ j= 1,\dots,A$) and transform
back to the image domain.
\begin{itemize}
\item Apply a threshold in the image domain.
\item Add each point $x$ above the threshold to the surfel set
  as the pair $(x,\kt)$.
\end{itemize}
\item Given the set of all surfels, use the algorithm of
\cite{meleslie:curveppaer} to reconstruct the set of curves
$\gamma_j(t)$ that define the discontinuities.
\end{enumerate}

The last step is described briefly in section
\ref{sec:connectingsurfels}, while the bulk of the paper is devoted to the
surfel extraction procedure itself.

\section{Large $k$ Asymptotics of  $\fstep{j}(\vk)$}
\label{sec:largeKAsymptotics}

In this section we wish to compute the large $k$ asymptotics of $\fstep{j}(\vk)$. In particular, we will show that $\fstep{j}(\vk)$ is dominated solely by the parts of $\gamma_{j}(t)$ where $\gamma_{j}'(t) \cdot \vec{k} = 0$, i.e. where $k \perp T_{j}(t)$.

We can use Green's theorem to rewrite the Fourier transform of $\step{j}(x)$ as follows. Let $F(\vk,x) = (F_1(\vk,x),F_2(\vk,x)) = -i \abs{\vk}^{-2} e^{i \vk \cdot x} \kp$ with $\kp=[-k_{2},k_{1}]^{T}$. Then by Green's theorem:
\begin{multline}
  \label{eq:exactTransformOfStepFunction}
  \fstep{j}(\vk) = \int \int_{\Omega_{j}} e^{i \vk \cdot x} dx_{1} dx_{2} =
  \int \int_{\Omega_{j}} \partial_{x_{1}} F_{2}(\vk,x) - \partial_{x_{2}} F_{1}(\vk,x) dx_{1} dx_{2}\\
  = \int_{\Sone} F(\vk, \gamma_{j}(t)) \cdot \frac{d \gamma_{j}(t)}{dt} dt =
  \frac{1}{ i\absSmall{\vk}^{2}}\int_{\Sone} e^{i \vk \cdot \gamma_{j}(t)} \kp \cdot \gamma_{j}'(t) dt
\end{multline}
(with $\Omega_{j}$ the region bounded by $\gamma_{j}$). This trick is taken from \cite{201624}.

We can now use stationary phase to approximate $\fstep{j}(\vk)$ for large $k$. For this,
express $k$ in polar coordinates $(\kr,\kt)$, fix a direction $\kt$, and consider what happens
as $\kr$ becomes large. As we remarked earlier, the phase $\vk \cdot \gamma_{j}(t)$ becomes stationary only when $\vk \cdot \gamma_{j}'(t)=0$ or $\vk \cdot T_{j}(t) =0$. This is precisely where $k$ points normal to the curve, and it is these locations that dominate $\fstep{j}(\vk)$:

\begin{proposition}
  \label{prop:LargeKAsymptoticsOfImage}
  Let $t_{j}(\vk)$ correspond to the value of $t$ at which $N_{j}(t) \parallel k$ and $N_{j}(t) \cdot \vec{k} > 0$ (i.e. the normal to $\gamma_{j}(t)$ points in the direction $k$). Then:
  \begin{multline}
    \label{eq:asymptoticExpansionOfStepFunction}
    \sum_{j=0}^{M-1} \ImgCoeff{j} \fstep{j}(\vk)
    = \sum_{j=0}^{M-1} \ImgCoeff{j} \left[
      \frac{ e^{i \vk \cdot \gamma(t_{j}(\vk))}}{\abs{\vk}^{3/2}} \frac{\sqrt{\pi} }{\sqrt{\kappa_{j}(t_{j}(\vk))} } + \frac{ e^{i \vk \cdot \gamma(t_{j}(-k))}}{\abs{\vk}^{3/2}} \frac{\sqrt{\pi} }{ \sqrt{\kappa_{j}(t_{j}(-\vk))} }
      \right]\\
      + \frac{E(\vk)}{\abs{\vk}^{2}}
  \end{multline}
  where $E(\vk) \leq \ConstGeo$
%%\begin{equation}
%    \label{eq:asymptoticExpansionOfStepFunctionRemainder}
%   \abs{E(\vk)} \leq  A_{geo}
% \end{equation}
  with
 \begin{equation}
    \label{eq:asymptoticExpansionOfStepFunctionRemainder2}
    \ConstGeo = M \contrastMax\left( 4 + \frac{8 \curvemax}{ \pi \curvemin} + 2 \sqrt{ \frac{2 \curvemax}{ \curvemin} } \right) + 3 \contrastMax\frac{\sup_{j} \norm{\gamma_{j}'''(t)}{L^{\infty}} }{ \curvemin} \sum_{j} \mathrm{arclength}(\gamma_{j}) \, .
  \end{equation}
  \end{proposition}

Note that $\ConstGeo$ incorporates both geometric and contrast information about the image itself.
This result is proved carefully in Appendix \ref{sec:stationaryPhaseCalculations}. The basic idea behind the proof is simple, however. Set $k=\kr \ktv$, with $\ktv$ fixed. Then:
\begin{equation*}
  \int_{\Sone} e^{i \vk \cdot \gamma_{j}(t)} \kp \cdot \gamma_{j}'(t) dt = \int_{\Sone} e^{i \kr f(t)} \kp \cdot \gamma_{j}'(t) dt
\end{equation*}
with $f(t) = \ktv \cdot \gamma_{j}(t)$. The phase function $f(t)$ is stationary when $\ktv \cdot \gamma_{j}'(t) = 0$, or equivalently the place where $\ktv \cdot N_{j}(t) = \pm 1$ (i.e. $t_{j}(\vk)=t_{j}(\ktv)$). By Assumption \ref{ass:curvatureBoundedBelow}, we find that $\gamma_{j}''(t) = \kappa_{j}(t) N_{j}(t)$ is nonzero. We restrict the consideration to an interval $I = [t_{j}(\ktv) - \epsilon, [t_{j}(\ktv) + \epsilon]$ and apply stationary phase:
\begin{multline*}
  \int_{I} e^{i \kr f(t)} \kp \cdot \gamma_{j}'(t) dt \\
  \sim
  \sqrt{\frac{\pi}{\kr f''(t_{j}(\ktv))}}e^{i \kr f(t_{j}(\ktv))} \kp \cdot \gamma_{j}'(t_{j}(\ktv)) + \textrm{remainder} \\
  = \kr^{1/2} \sqrt{\frac{\pi}{\kappa_{j}(t_{j}(\ktv))}} e^{i \vk \cdot \gamma_{j}(t_{j}(\vk))} + \textrm{remainder}
\end{multline*}
Proving Proposition \ref{prop:LargeKAsymptoticsOfImage} is done by adding this result up over all the curves, and all points of stationary phase, and estimating the remainder.

\subsection{What if $\curvemin = 0$?}

It is important to note that all is not lost when $\curvemin = 0$. In this case, although the $\emph{coefficient}$ on $k^{-3/2}$ in \eqref{eq:asymptoticExpansionOfStepFunction} becomes singular, the asymptotics of $\fstep{j}(\vk)$ do not.

In this case, what happens is that the leading order behavior becomes $O(\kr^{-1-1/n})$, where $n$ is the order of the first non-vanishing derivative. This can be seen relatively easily from stationary phase, although rigorous justification requires a long calculation. However, the limiting case $n=\infty$ (a straight line) is easy enough to treat.

\begin{proposition}
  Let $\gamma(t)=\va + (\vb - \va)t$ for $t \in [0,1]$. Then:

  \begin{multline}
    \label{eq:straightLineSegment}
    \frac{1}{i \abs{\vk}^{2} } \int_{0}^{1} e^{i \vk \cdot \gamma(t)} \kp \cdot \gamma'(t) dt \\
    =
    \left\{
    \begin{array}{ll}
    \frac{-1}{\abs{\vk}^{2}} \frac{ \kp \cdot (\vb - \va)}{ \vk \cdot{ (\vb - \va)} }\left[
        e^{i \vk \cdot \vb} - e^{i \vk \cdot \va}
      \right], & \vk \cdot (\vb - \va) \neq 0 \\
      \frac{e^{i \vk \cdot \va } }{\abs{\vk}^{2}}\kp \cdot (\vb - \va), &k \cdot (\vb - \va) = 0
    \end{array}
    \right. \\
    =
    \left\{
    \begin{array}{ll}
      O(\kr^{-2}), & k \cdot (\vb - \va) \neq 0 \\
      e^{i k \cdot \va } O(\kr^{-1}), &k \cdot (\vb - \va) = 0
    \end{array}
    \right.
  \end{multline}
\end{proposition}

\begin{proof}
  If $\vk \cdot (\vb - \va) \neq 0$, then:
  \begin{multline*}
    \int_{0}^{1} e^{i k \cdot \gamma(t)}  \kp \cdot \gamma'(t) dt = \int_{0}^{1} e^{i k \cdot ((\vb - \va) t + \va)} \kp \cdot (\vb - \va) dt \\
    = \kp \cdot (\vb - \va) e^{i k \cdot \va} \int_{0}^{1} e^{i k \cdot ((\vb - \va) t)} dt = \kp \cdot (\vb - \va) e^{i k \cdot \va} \left[ \frac{e^{i k \cdot (\vb - \va)}}{ i k \cdot (\vb - \va)} - \frac{1}{ i k \cdot (\vb - \va)}\right] \\
    = \frac{ \kp \cdot (\vb - \va) }{i k \cdot (\vb - \va)} \left[ e^{i k \cdot \vb } - e^{i k \cdot \va}\right]
  \end{multline*}
  Multiplying by $(i \abs{\vk}^{2})^{-1}$ yields the result we seek.

  The asymptotics are straightforward to compute (from the second line of \eqref{eq:straightLineSegment}), but note that the constant in the $O(\kr^{-2})$ term in \eqref{eq:straightLineSegment} is not uniform in $k \angle (\vb - \va)$.
\end{proof}

\section{Directional Filters}
\label{sec:directionalFilters}

We are now in a position to build the filter operators $\dfilter{\theta}{\alpha}$
of \eqref{eq:filterDef}, which will allow us to extract
edge information from the signal. We demand that the radial filter takes the form
\begin{subequations}
  \begin{equation}
    \label{eq:formOfRadialFilter}
    \RadFilt(\kr) = \RadFiltCentered(\kr - (\kmax + \ktex)/2) + \RadFiltCentered(\kr + (\kmax + \ktex)/2)
  \end{equation}
  where $\RadFiltCenteredX(r)$ is a positive, symmetric function. This means that the
  $1$-dimensional inverse Fourier transform of $\RadFilt(\kr)$ is
  \begin{equation}
    \label{eq:inverseFourierTransformOfRadialFilter}
    \check{\RadFilt}(r) = \cos\left(\frac{\kmax + \ktex}{2}r \right) \RadFiltCenteredX(r)
  \end{equation}
\end{subequations}
The intuition behind this operator is the following. Multiplying by $\RadFilt(\kr) \AngFilt(\kt - \theta)$ localizes on the region $k \approx (\kr \cos \theta, \kr \sin \theta)$. Dropping all but one of the terms in Proposition \ref{prop:LargeKAsymptoticsOfImage}, we find that:
\begin{equation*}
  \RadFilt(\kr) \AngFilt(\kt - \theta) \fstep{k}(\vk) \approx e^{i k \cdot \gamma_{j}(t_{j}(\vk))} \kappa_{j}(t_{j}(\vk))^{-1/2} A(\vk) \, .
\end{equation*}
The term $A(\vk)$ incorporates both the $1/|\vk|^{3/2}$ decay from the asymptotics of $\fstep{k}(\vk)$ and the localization to $\kt$ and large $\kr$ from the filters.

Therefore, if we inverse Fourier transform, we will obtain
\begin{equation*}
  \ifourier [e^{i k \cdot \gamma_{j}(t_{j}(\vk))} \kappa_{j}(t_{j}(\vk))^{-1/2} A(\vk)] \approx \kappa_{j}(t_{j}(\vk))^{-1/2} \check{A}(x - \gamma_{j}(t_{j}(\vk)))
\end{equation*}
Provided $\check{A}(x)$ is a sharply localized bump function, this will be a bump located at $\gamma_{j}(t_{j}(\vk))$. Of course, this calculation is not exactly correct, and is presented merely to obtain intuition. We will go through the details shortly, but require a few definitions.

\begin{definition}
  Given $\RadFilt(\kr)$, define the auxiliary functions:
  \begin{subequations}
    \begin{eqnarray}
      \check{\RadFilt}(r) & = & \int e^{-i r \kr} \RadFilt(\kr) d\kr \label{eq:defOfcheckf}\\
      \RadFiltXSup(R)&=&\sup_{r > R} \abs{\check{\RadFilt}(r)} \label{eq:defOfBarf}
    \end{eqnarray}
  \end{subequations}
\end{definition}

\begin{definition}
  \label{def:arcsOfSpecificAngle}
  Define the set $\arcs{\theta}{\alpha} \subset \RR^{2}$ to be the set of points where some $\gamma_{j}(t)$ has normal pointing in the direction $[\theta-\alpha,\theta+\alpha]$,
 i.e.:
  \begin{subequations}
    \begin{equation}
      \label{eq:arcsDef}
      \arcs{\theta}{\alpha} = \bigcup_{j=0}^{M-1} \left[ \{ \gamma_{j}(t) \}_{t \in [t_{j}(\theta-\alpha),t_{j}(\theta+\alpha)]} \cup \{ \gamma_{j}(t) \}_{t \in [t_{j}(\theta-\alpha + \pi),t_{j}(\theta+\alpha + \pi)]} \right]
    \end{equation}
    We also define $\arcsMinusOne{\theta}{\alpha}{j}$ to be the set of arcs excluding $\gamma_{j}(t_{j}(\theta))$.
    \begin{equation}
      \label{eq:arcsMinusOneDef}
      \arcsMinusOne{\theta}{\alpha}{j} = \bigcup_{i \neq j} \left[ \{ \gamma_{i}(t) \}_{t \in [t_{i}(\theta-\alpha),t_{i}(\theta+\alpha)]} \cup \{ \gamma_{i}(t) \}_{t \in [t_{i}(\theta-\alpha + \pi),t_{i}(\theta+\alpha + \pi)]} \right]
    \end{equation}
  \end{subequations}
\end{definition}
The goal of our directional filters is to approximate the location of $\arcs{\theta}{\alpha}$. That is to say we want $\dfilter{\theta}{\alpha} \Img$ to be large near $\arcs{\theta}{\alpha}$ and small away from it. The decay in the tangential direction away from a point in $\arcs{\theta}{\alpha}$
 is at least of the order $O(r^{-1})$. To ensure that the decay in the normal direction is as fast,
 we consider the case when $\RadFiltX(r) = O(r^{-1})$ (see Fig.
\ref{fig:directionalFilterSimpleUse} \, d).

 We then have the following result which proves the directional filters approximate the location of $\arcs{\theta}{\alpha}$.
\begin{theorem}
  \label{thm:DirectionalFilterYieldsEdges}
  \begin{subequations}
    Suppose that $\alpha$ satisfies
    \begin{equation}
      \label{eq:geometricConstraintMaxOfFilter}
      \frac{\alpha^{2} \left(\kmax + \ktex \right)}{ \curvemin} \leq \pi
    \end{equation}
    and $\RadFilt(\kr)$ satisfies \eqref{eq:formOfRadialFilter} as well as the following decay condition:
    \begin{equation}
      \label{eq:decayOfRadFilterInX}
      \RadFiltXSup(r) \leq \frac{\ConstRadFiltDecay}{r}
    \end{equation}

    Then away from $\arcs{\theta}{\alpha}$, the directional filter has the following decay:
    \begin{multline}
      \label{eq:DirectionalFilterYieldsEdges}
      \abs {\dfilter{\theta}{\alpha} \left[\sum_{j=0}^{M-1} \ImgCoeff{j} \fstep{j}(x) \right] } \leq
      \frac{\ConstTangentNormal}{d(x, \arcs{\theta}{\alpha} ) }\\
      + \norm{\kr^{-1/2} \RadFilt(\kr)}{L^{1}(\RR, d\kr)} \norm{\AngFilt(\kt)}{L^{1}(\Sone,d\kt)} \, \ConstGeo
    \end{multline}
    where
    \begin{multline}
      \label{eq:defOfConstTangentNormal}
      \ConstTangentNormal = \frac{2M \contrastMax \sqrt{2\pi}}{\sqrt{\curvemin}
        \cos(2\alpha)} \max\Bigg\{
      \ConstRadFiltDecay \norm{\AngFilt(\kt)}{L^{1}(\Sone,d\kt)},\\
      \Bigg( 2 \norm{\RadFilt(\kr)/\kr}{L^{1}(\RR, d\kr)}
      \left[\norm{\AngFilt'(\kt)}{L^{1}}
        + \norm{\AngFilt(\kt)}{L^{1}}  \frac{\norm{\gamma_{j}'''(t)}{L^{\infty}}}{2 \curvemin^{2}} \right]\\
      + \sqrt{\curvemin} \norm{\check{\RadFilt}(z)(z+1)}{L^{\infty}}
      \Bigg) \Bigg\}
    \end{multline}
    At the point $x = \gamma_{j}(t_{j}(\theta))$
    \begin{equation}
      \label{eq:DirectionalFilterLargeAtEdge}
      [\theta \cdot N_{j}(t)] \dfilter{\theta}{\alpha} \left[\sum_{j=0}^{M-1} \ImgCoeff{j} \fstep{j}(x) \right]
      \geq \SurfelThreshold
    \end{equation}
    with
    \begin{multline}
        \label{eq:SurfelThresholdDef}
        \SurfelThreshold
        \equiv \sqrt{\frac{\pi}{2 \curvemax}} \contrastMin \inf_{r \in [-\alpha^{2}/2\curvemin, \alpha^{2}/2\curvemin]} \RadFiltCentered(r)\\
      - \left[ \frac{\ConstTangentNormal (2M-1) }{2M \, \curvesep}
      + \norm{\kr^{-1/2} \RadFilt(\kr)}{L^{1}(\RR, d\kr)} \right]  \, \ConstGeo
      \end{multline}
  \end{subequations}
\end{theorem}

We postpone the proof of this result to Appendix \ref{sec:DirectionalFiltersLeadingOrderAsymptotocs}, where we compute the leading order asymptotics of a single segment of the curve and put the various pieces together.

\begin{remark} \label{rmk:competition}
While the expressions above
are somewhat involved, for the filters described in the next section, \eqref{eq:DirectionalFilterYieldsEdges} simplifies to
 \begin{multline}
      \label{eq:DirectionalFilterYieldsEdgesSimple}
      \abs {\dfilter{\theta}{\alpha} \left[\sum_{j=0}^{M-1} \ImgCoeff{j} \fstep{j}(x) \right] } \leq
      \frac{\ConstTangentNormal}{d(x, \arcs{\theta}{\alpha} ) } +
      O \Bigg(\frac{1}{\kmax^{1/2} + \ktex^{1/2}} \Bigg).
    \end{multline}
Note that the second term in this estimate becomes negligible as $\kmax$ increases.
The term  $\ConstTangentNormal$ determines the rate of decay of the directional filter away
from a surfel and it should be as small as possible (to minimize noise), for which
the norm of  $\RadFilt$ should be small. At the same time, we want
 $\SurfelThreshold$ to be as large as possible (to maximize signal).
 For this, the norm of  $\RadFilt$ should be large. We will need to balance this competition.
\end{remark}

\begin{remark}
  If the decay of $\RadFiltXSup(r)$ is faster than \eqref{eq:decayOfRadFilterInX}, we can of course get sharper results. However, this would require extra geometric conditions and would require Theorem \ref{thm:DirectionalFilterYieldsEdges} to make distinctions concerning the direction from $x$ to $\arcs{\theta}{\alpha}$.  To simplify the exposition of this paper, such considerations will be reported at a later date.
\end{remark}

\subsection{The estimate is suboptimal}
\label{sec:suboptimalEstimateExample}

\begin{figure}
\setlength{\unitlength}{0.240900pt}
\ifx\plotpoint\undefined\newsavebox{\plotpoint}\fi
\sbox{\plotpoint}{\rule[-0.200pt]{0.400pt}{0.400pt}}%
\includegraphics[scale=0.75, bb=100pt 175pt 200pt 450pt]{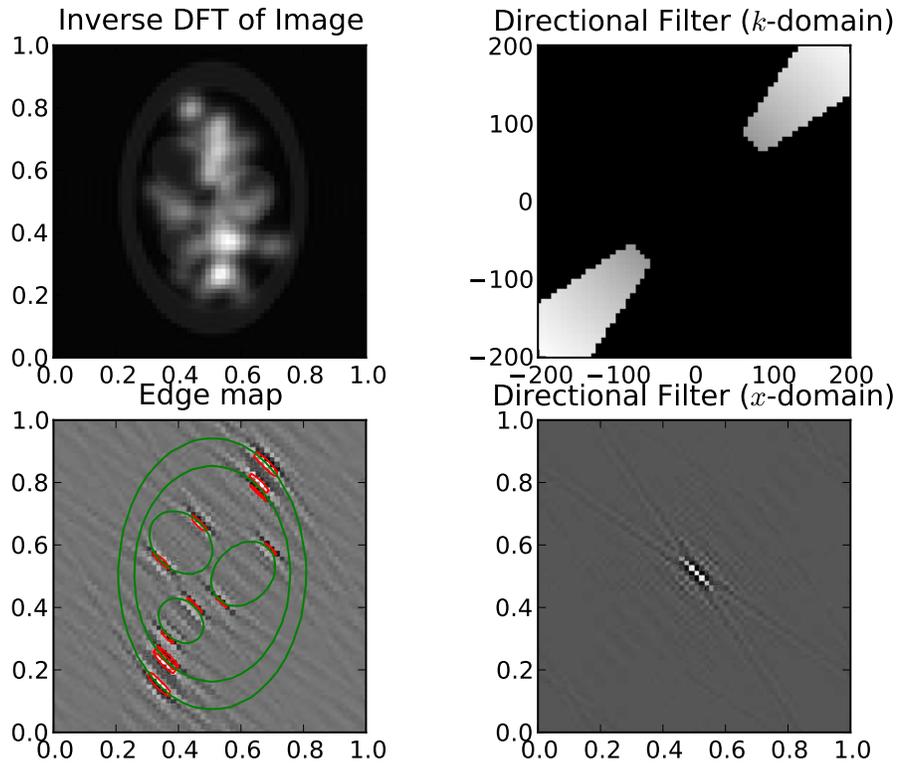}

\caption{An illustration of the sharp directional filter. The first panel shows the image. The second panel shows the directional filter in $k$-space. The third panel shows the edge map $[\dfilter{\theta}{\alpha} \Imgk](x)$ with $\theta=\pi/4$ and $\alpha=\pi/16$. The red lines are the $\abs{[\dfilter{\theta}{\alpha} \Imgk](x)}=2.4$ contour lines, while the green lines are the actual (analytically known) edges of the image. The fourth panel shows the directional filter in the $x$-domain.}
\label{fig:directionalFilterSimpleUse}
\end{figure}

\begin{figure}
\setlength{\unitlength}{0.240900pt}
\ifx\plotpoint\undefined\newsavebox{\plotpoint}\fi
\sbox{\plotpoint}{\rule[-0.200pt]{0.400pt}{0.400pt}}%
\subfigure{
\includegraphics[scale=0.65, bb=100pt 175pt 200pt 500pt]{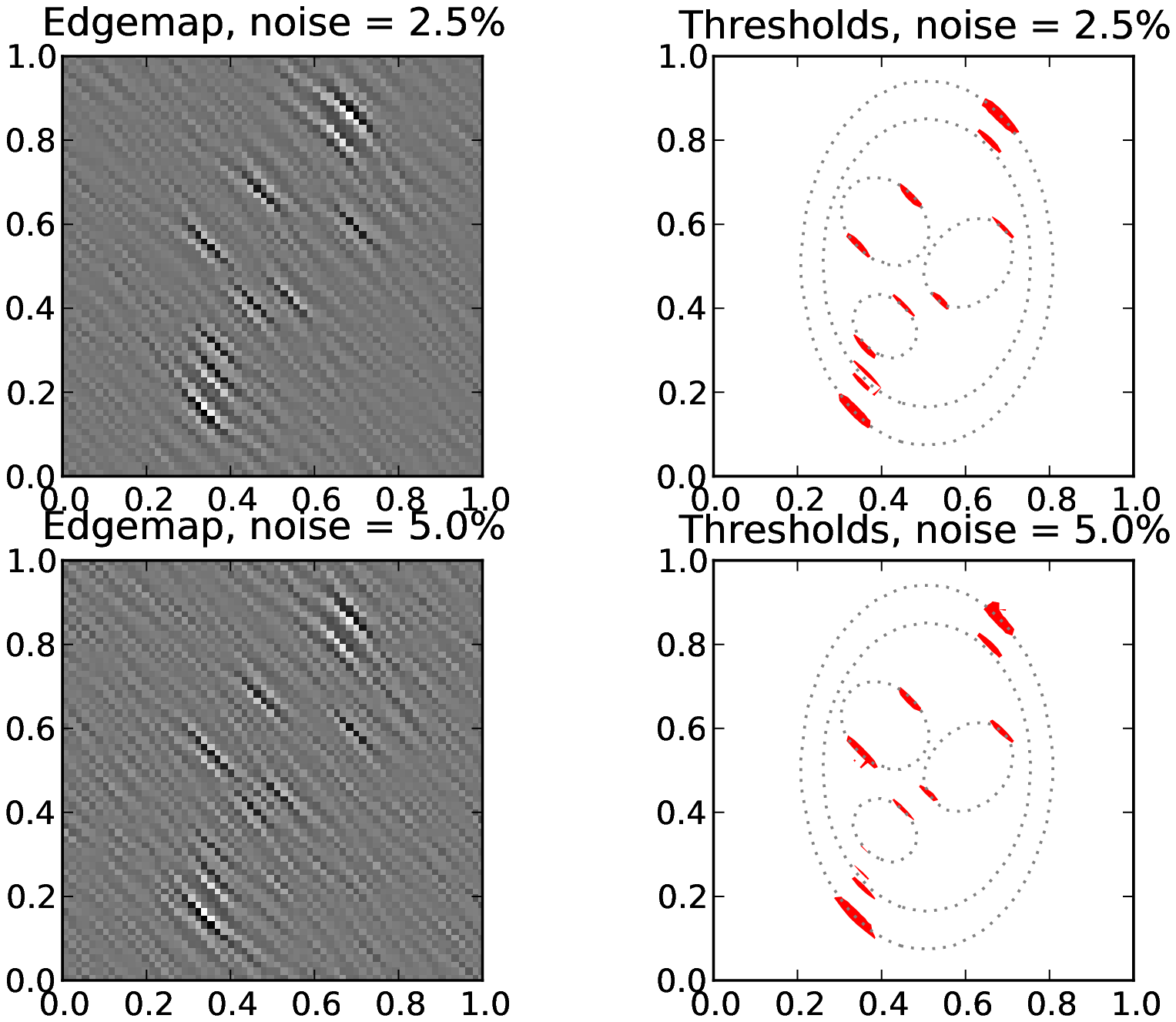}
}
\newline
\subfigure{
\includegraphics[scale=0.65, bb=100pt 175pt 200pt 565pt]{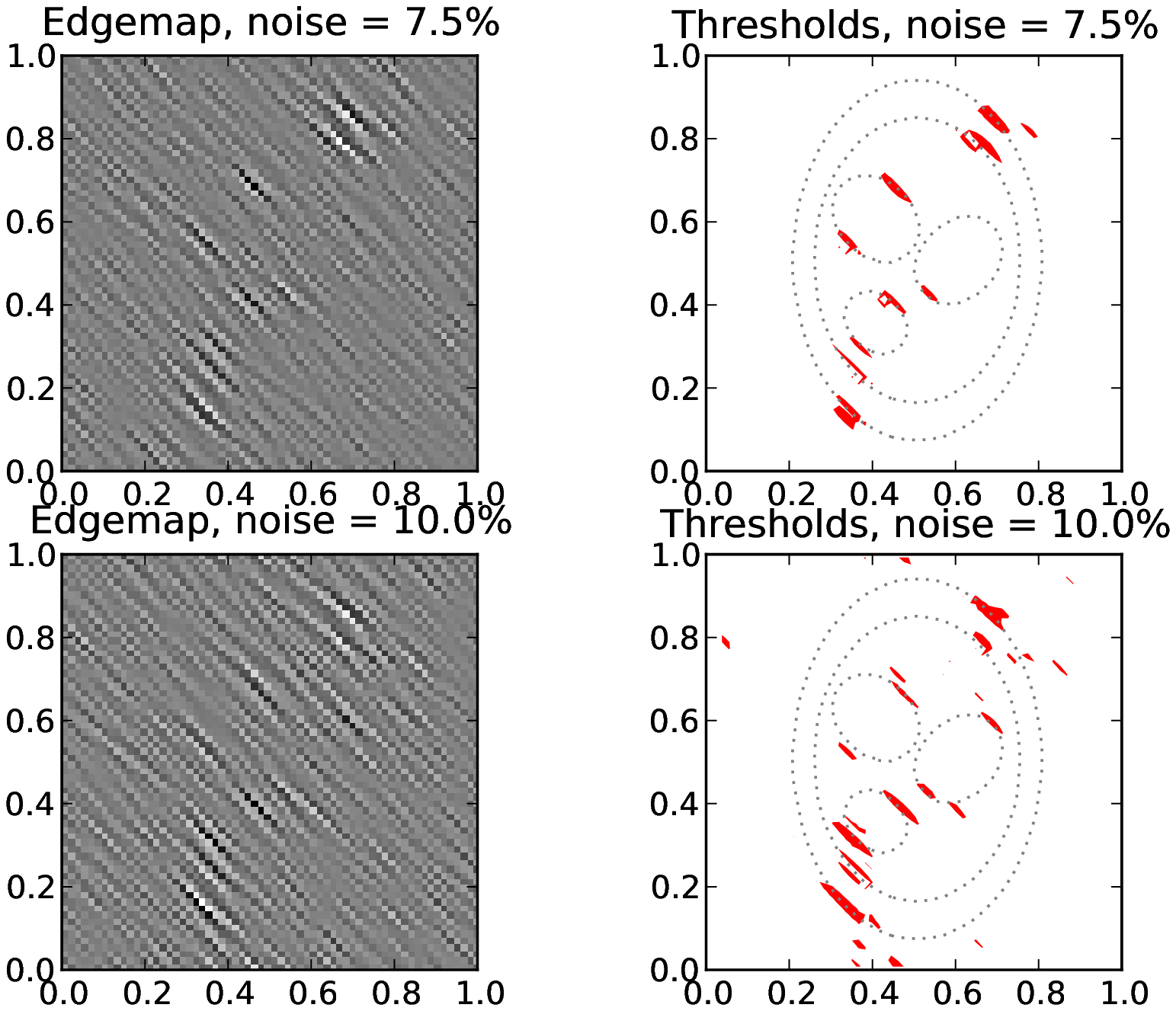}
}

\caption{An illustration of the directional filters in the presence of noise. The image and parameters are the same as in Figure \ref{fig:directionalFilterSimpleUse}. The figures in the left column are plots of $[\dfilter{\theta}{\alpha} \Imgk](x)$. The gray dotted lines in the right column illustrate the discontinuities of the image (which we know exactly, c.f. Section \ref{sec:PhantomExplanation}), while the red regions are the places where $\abs{[\dfilter{\theta}{\alpha} \Imgk](x)}=2.4$. Noise amplitude is described as a percentage of the total image strength. Noise begins to corrupt the data when it reaches 7.5\%, and becomes a serious problem at 10\%. }
\label{fig:directionalFilterSimpleUseNoisy}
\end{figure}

Let us consider the application of our directional filters on a simple numerical example. We consider a $64 \times 64$ pixel image on $[0,1]^{2}$. The image is taken to be the phantom described in Section \ref{sec:PhantomExplanation}, and the parameters are $\kmax = 64 \pi \approx 201.1$, $\ktex = 32 \pi \approx 100.5$, $\alpha=\pi/16$, $\curvemin=0.1$ and $\norm{\gamma_{j}'''(t)}{L^{\infty}} \leq 5$. With this set of parameters, $\ConstTangentNormal = 0.3$ if we use the windows $\AngFilt(\kt) = (2\alpha)^{-1} 1_{[\theta-\alpha, \theta+\alpha]}(\kt)$ and $\RadFilt(\kr) = (\kmax-\ktex)^{-1} 1_{[\ktex, \kmax]}(\kr)$.
The first term on the right-hand side in  \eqref{eq:DirectionalFilterYieldsEdges} is, therefore, approximately
$0.3/(1/64) \approx 25$ if we move one pixel away from a surfel. The second term, however, is approximately $85$ so that we have no reason to expect that the directional filters will yield any useful information. However, numerical experiments show that they do in fact work even in this case. Figure \ref{fig:directionalFilterSimpleUse} shows that the directional filters do yield the location of the edges. Figure \ref{fig:directionalFilterSimpleUseNoisy} shows that even in the presence of noise (up to $7.5\%$ of the total image energy), the directional filters still yield correct results.

\begin{figure}
\setlength{\unitlength}{0.240900pt}
\ifx\plotpoint\undefined\newsavebox{\plotpoint}\fi
\sbox{\plotpoint}{\rule[-0.200pt]{0.400pt}{0.400pt}}%
\includegraphics[scale=0.75, bb=100pt 200pt 200pt 520pt]{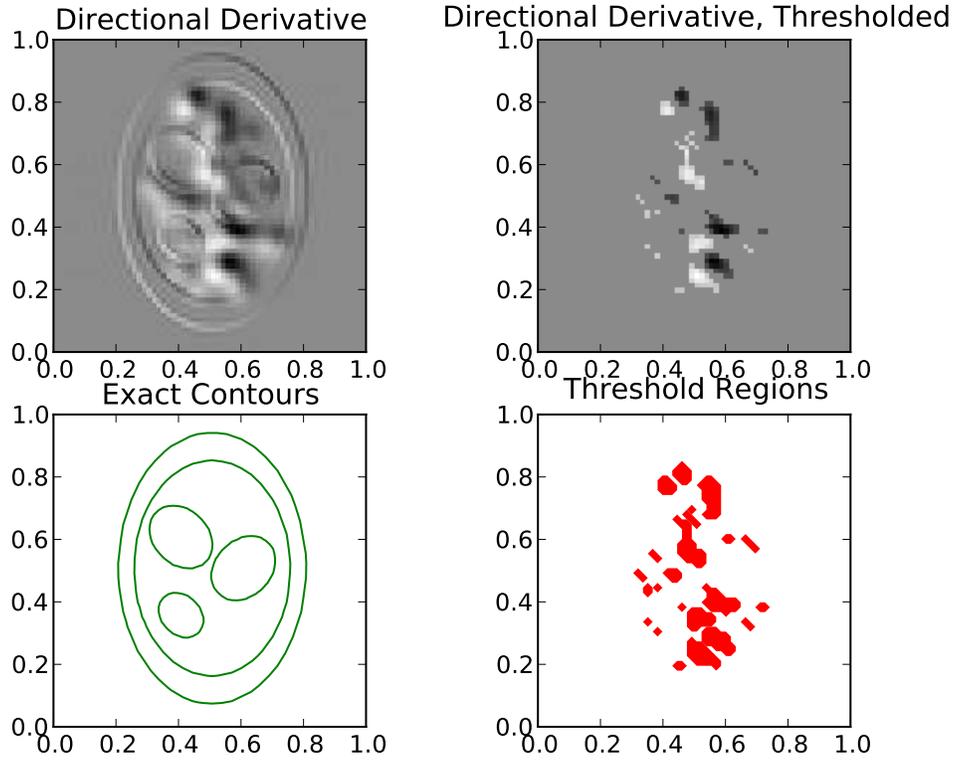}

\caption{Result of edge detection based on measuring the directional derivative, i.e. $\partial_{\pi/4} \Img(x)$. The green lines represent the analytically known edges of the image, while the red regions are the regions where $\abs{\partial_{\pi/4} \Img(x)} \geq 150$.}
\label{fig:ComparisonToDerivativeDetector}
\end{figure}

It is also useful to compare the directional filters to a naive algorithm, namely computing the directional derivative. As is apparent from Figure \ref{fig:ComparisonToDerivativeDetector}, this naive method does not perform as well as our algorithm.

%%%\newpage

%%\newpage

%%%\newpage
\begin{figure}
\setlength{\unitlength}{0.240900pt}
\ifx\plotpoint\undefined\newsavebox{\plotpoint}\fi
\sbox{\plotpoint}{\rule[-0.200pt]{0.400pt}{0.400pt}}%
\subfigure{
\includegraphics[scale=0.65, bb=100pt 175pt 200pt 500pt]{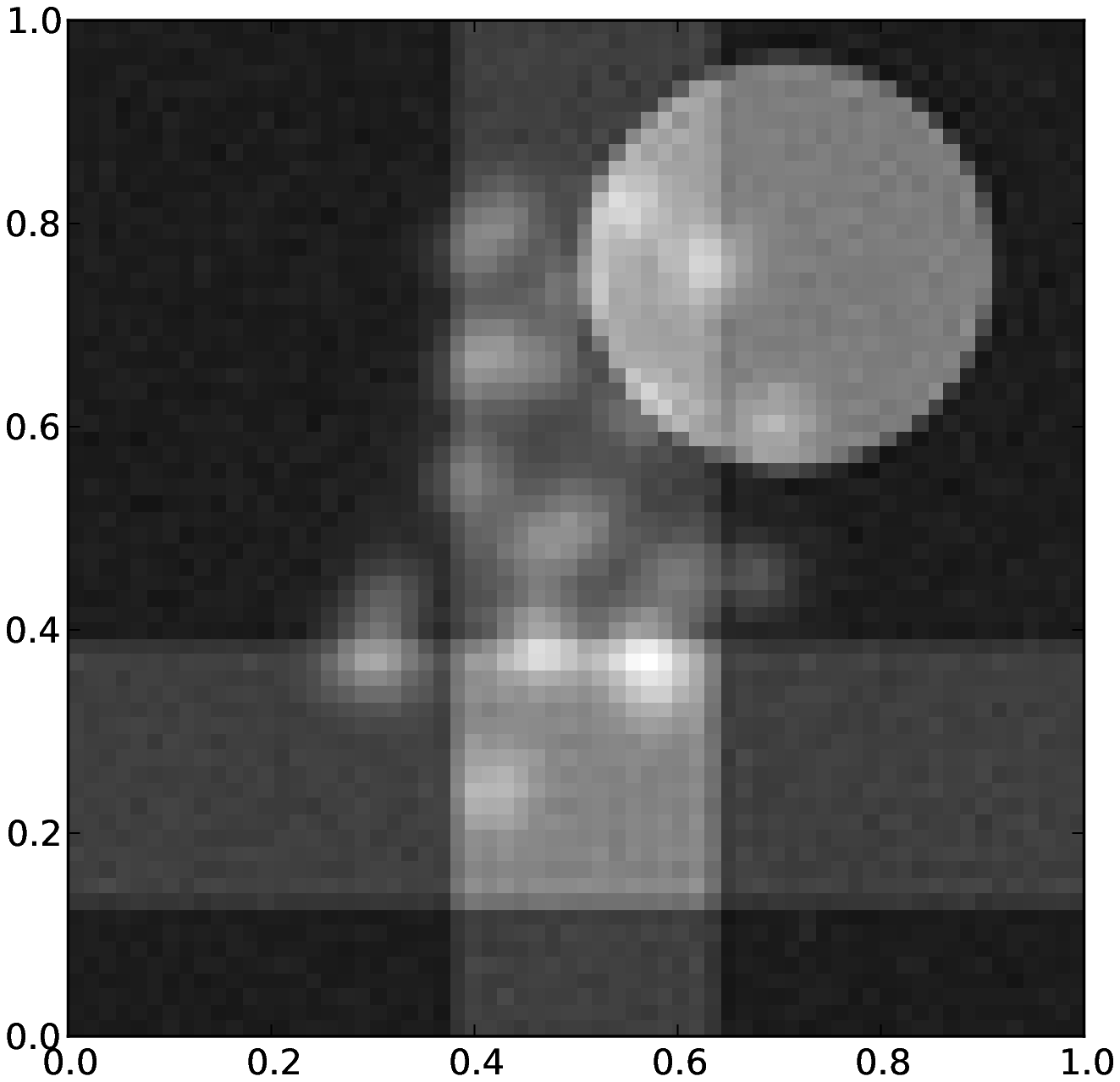}
}
\newline
\subfigure{
\includegraphics[scale=0.65, bb=100pt 175pt 200pt 565pt]{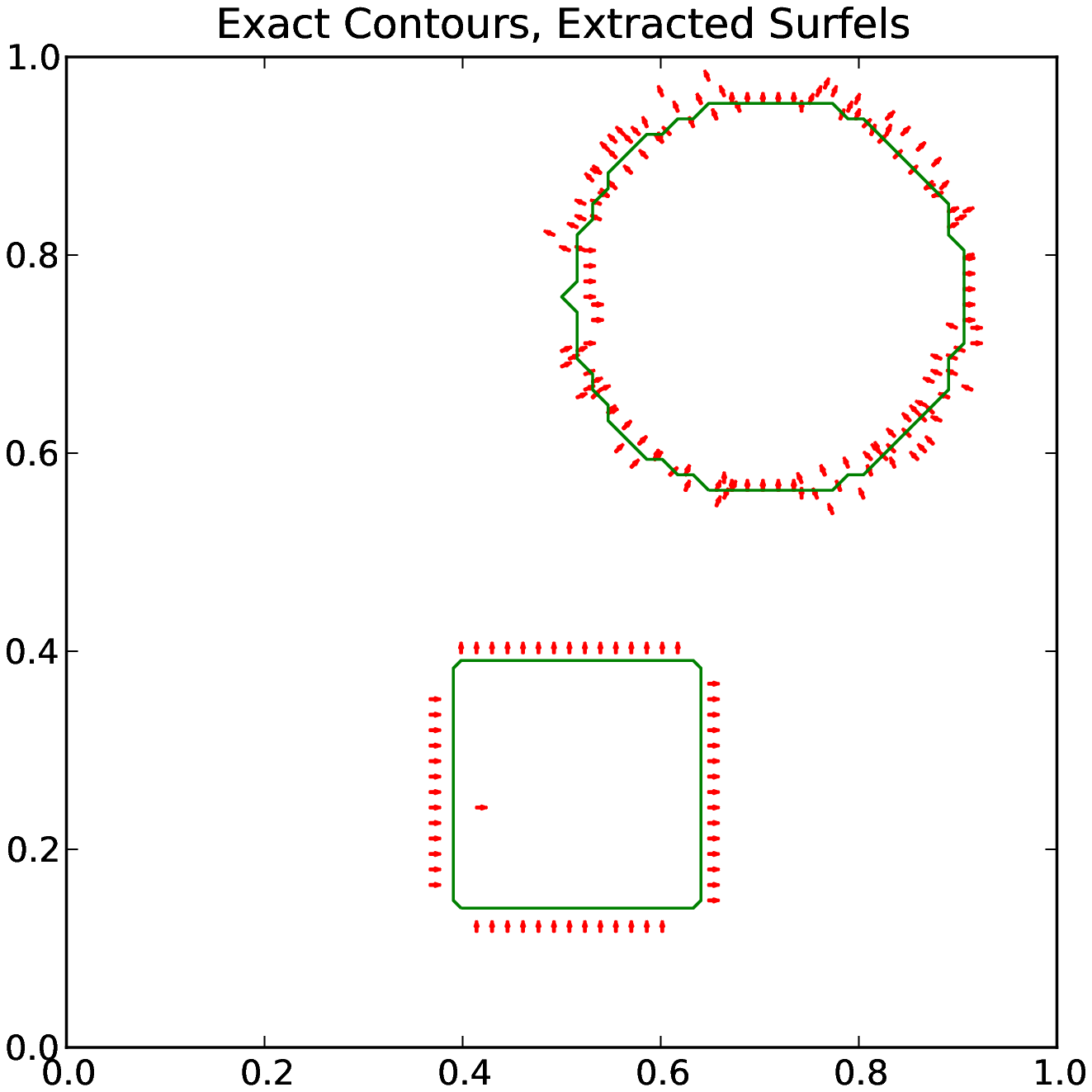}
}

\caption{Result of edge detection applied to an image where Assumption \ref{ass:curvatureBoundedBelow} is violated. The green lines in the second panel indicate the actual discontinuities, while the red arrows indicate extracted surfels (at angles $\theta=0$ and $\theta=\pi/4$. The vertical and horizontal Gray bars are artifacts from using the DFT to reconstruct a discontinuous image.}
\label{fig:squareExample}
\end{figure}

\subsection{Parabolic Scaling}

The choice of $\alpha$ is an important one. We want $\alpha$ to be as small as possible, since this will give us better angular resolution.

On the other hand, the constant bounding the size of the filtered image away from the edges is directly proportional to $\alpha^{-1}$, as we will show shortly. For simplicity, let us impose the constraint that
\begin{equation}
  \label{eq:angularFilterNormIs1}
  \norm{\AngFilt(\kt)}{L^{1}(\kt, d\kt)} = 1.
\end{equation}
Now let us take $\AngFilt(\kt)$ to be a fixed filter scaled with $\alpha$, i.e. $\AngFilt(\kt) = \alpha^{-1} V(\kt/\alpha)$, where $V(\kt)$ is supported in $[-1,1]$ and $\norm{V(\kt)}{L^{1}(d\kt)} = 1$. With this choice of $\AngFilt(\kt)$, we find that:
\begin{equation}
  \label{eq:angularFilterDerivNormIsAlpha}
  \norm{\AngFilt'(\kt)}{L^{1}(\Sone,d\kt)} = \norm{\alpha^{-1} V'(\kt/\alpha)}{L^{1}(\Sone, d\kt)} = \alpha^{-1} \norm{V'(\kt/\alpha)}{L^{1}(\Sone, d\kt)}
\end{equation}

Substituting this into \eqref{eq:defOfConstTangentNormal} and assuming $\alpha$ to be very small yields:
\begin{multline}
  \ConstTangentNormal = \frac{\sqrt{2\pi}}{\sqrt{\curvemin} \cos(2\alpha)}
  \Bigg(
  2 \norm{\RadFilt(\kr)/\kr}{L^{1}(\RR, d\kr)} \\
  \times \left[\alpha^{-1} \norm{V'(\kt/\alpha)}{L^{1}(\Sone, d\kt)}
    + \norm{\AngFilt(\kt)}{L^{1}}  \frac{\norm{\gamma_{j}'''(t)}{L^{\infty}}}{2 \curvemin^{2}} \right]\\
  + \sqrt{\curvemin} \norm{\check{\RadFilt}(z)(z+1)}{L^{\infty}}
  \Bigg) = O\left(\curvemin^{-1/2} \alpha^{-1} \norm{\RadFilt(\kr)/\kr}{L^{1}(\RR, d\kr)} \right)
\end{multline}
To prevent $\ConstTangentNormal$ from being too large, we want to ensure that $\alpha$ is as large as possible.

However, to ensure that the filter is sufficiently large \emph{on the edges}, we need to prevent $\alpha$ from being too large (Theorem \ref{thm:DirectionalFilterYieldsEdges}).
In particular, we require (c.f. \eqref{eq:geometricConstraintMaxOfFilter})
\begin{equation}
  \label{eq:19}
  \alpha \leq \sqrt{\frac{\pi \curvemin}{\kmax + \ktex}}
\end{equation}
To satisfy both these constraints, we simply replace the ``$\leq$'' in \eqref{eq:19} by ``$=$''. This yields the standard parabolic scaling used elsewhere
\cite{MR2165380,MR2012649,fefferman:parabolicScaling,seggerSoggeStein:parabolicScaling,stein:harmonicAnalysis}.
\begin{remark}
  The requirement \eqref{eq:geometricConstraintMaxOfFilter} yields the standard parabolic scaling used to analyze line discontinuities in harmonic analysis   \cite{MR2165380,MR2012649,fefferman:parabolicScaling,seggerSoggeStein:parabolicScaling,stein:harmonicAnalysis}.
In the $x$-domain, the standard parabolic scaling uses elements with $\textrm{width} \sim \textrm{length}^{2}$. In the $k$-domain, this translates to $\textrm{width} \sim \sqrt{\textrm{length}}$.
\end{remark}

This also implies that for the filter to approach zero as $\ktex, \kmax \rightarrow \infty$, we require that $\norm{\RadFilt(\kr)/\kr}{L^{1}(\RR, d\kr)} = o(\ktex^{-1/2})$.

\subsection{Choosing the Filters} \label{sec:choosingfilters}

We now consider the simplest choices of window possible which satisfy our assumptions. We take $\RadFilt(\kr)$ to be a a step function, i.e.:
\begin{equation*}
  \RadFilt(\kr) = \frac{1}{2(\kmax - \ktex)} \left( 1_{[\ktex,\kmax]}(\kr) + 1_{[-\kmax,-\ktex]}(\kr) \right)
\end{equation*}
and we take $\AngFilt(\kt)$ to be a triangle function:
\begin{equation*}
  \AngFilt(\kt) = \frac{1}{\alpha^{2}}\max\{ \alpha - \abs{\kt}, 0 \}
\end{equation*}
In this case, it is straightforward to show:
\begin{eqnarray*}
  \norm{\kr^{-1/2} \RadFilt(\kr)}{L^{1}(\RR, d\kr)} & = & \frac{1}{\kmax^{1/2} + \ktex^{1/2}} \\
  \norm{\kr^{-1} \RadFilt(\kr)}{L^{1}(\RR, d\kr)} & = & \frac{\ln(\kmax/\ktex)}{\kmax - \ktex} \\
  \norm{ \RadFiltX(r)(r+1) }{L^{\infty}(\RR, dr)} & \leq & \frac{1}{\kmax - \ktex} \\
  \ConstRadFiltDecay & = &\frac{1}{\kmax - \ktex}\\
  \norm{\AngFilt(\kt)}{L^{1}(\Sone, d\kt)} & = & 1 \\
  \norm{\AngFilt'(\kt)}{L^{1}(\Sone, d\kt)} & = & 2\alpha^{-1} = 2 \sqrt{\frac{\kmax + \ktex}{\pi}}
\end{eqnarray*}
Then we have:
\begin{multline}
  \label{eq:stepWindowsC1}
  \ConstTangentNormal = \frac{\sqrt{2\pi}}{\sqrt{\curvemin} \cos(2\alpha) (\kmax - \ktex)}
  \max\Bigg\{1
  ,\\
  2 \ln(\kmax/\ktex) \Bigg( 2 \sqrt{\frac{\kmax + \ktex}{\pi \curvemin}} +
  \frac{\norm{\gamma_{j}'''(t)}{L^{\infty}}}{2 \curvemin^{2}}
  + \sqrt{\curvemin} \Bigg)
  \Bigg\}
\end{multline}
Similarly, the remainder term has the bound:
  \begin{multline*}
    \abs{\dfilter{\theta}{\alpha}\frac{E(\vk)}{\abs{k}^{2}}}
    \leq  \frac{1}{\kmax^{1/2} + \ktex^{1/2}} \norm{\AngFilt(\kt)}{L^{1}(\Sone,d\kt)} \\
    \times \ConstGeo \, ,
  \end{multline*}
  where $\ConstGeo$ is defined in (\ref{eq:asymptoticExpansionOfStepFunctionRemainder2}).

  Therefore, we find that:
  \begin{multline}
    \abs {\dfilter{\theta}{\alpha}
      \left[\sum_{j=0}^{M-1} \fstep{j}(x) \right] } \\
    = O\left(
      \frac{1}{d(x, \arcs{\theta}{\alpha} ) } \left[
        \frac{\ln(\kmax/\ktex) \sqrt{\kmax + \ktex}}{\kmax - \ktex}
        \right] + \frac{1}{\kmax^{1/2} + \ktex^{1/2}}
      \right)
  \end{multline}
  Thus, away from $\arcs{\theta}{\alpha}$, the directional filter approaches zero. The $O(\ln(\kmax/\ktex)/(\kmax - \ktex))$ factor is present because we are attempting to extract spatial information from a frequency band of width $\kmax - \ktex$. The $1/(\kmax^{1/2} + \ktex^{1/2})$ term is present because the asymptotic expansion we used to derive the directional filters decays only $\kr^{-1/2}$ faster than the leading order terms.

\section{Surfel Extraction}

Theorem \ref{thm:DirectionalFilterYieldsEdges} proves that, provided we choose the parameters correctly, directional filters will decay away from $\arcs{\theta}{\alpha}$. Extracting surfels from the filters is therefore simply a matter of choosing the parameters correctly and seeking local maxima.

We know that at the point $\gamma_{j}(t_{j}(\theta))$, \eqref{eq:DirectionalFilterLargeAtEdge} provides a lower bound on the size of the directionally filtered image. We also know that away from $\arcs{\theta}{\alpha}$, \eqref{eq:DirectionalFilterYieldsEdges} provides an upper bound on the size of the filtered image.

We wish to say that if the filtered image is ``large'', then we are near $\arcs{\theta}{\alpha}$, otherwise we are not. Therefore, we will need the lower bound at $\gamma_{j}(t_{j}(\theta))$ to be greater than the upper bound away from $\arcs{\theta}{\alpha}$:
\begin{multline}
  \label{eq:22}
  \sqrt{\frac{\pi}{2 \curvemax}} \contrastMin \norm{\AngFilt(\kt)}{L^{1}(d\kt)}
  \inf_{r \in [-\alpha^{2}/2\curvemin, \alpha^{2}/2\curvemin]}
  \RadFiltCentered(r)\\
  -  \left[ \frac{\ConstTangentNormal (2M-1) }{2M \, d(x, \arcsMinusOne{\theta}{\alpha}{j} ) }
    + \norm{\kr^{-1/2} \RadFilt(\kr)}{L^{1}(\RR, d\kr)} \norm{\AngFilt(\kt)}{L^{1}(\Sone,d\kt)}\right]
  \, \ConstGeo  \\
  \geq~
  \frac{\ConstTangentNormal}{d(x, \arcs{\theta}{\alpha} ) }\\
  +  \norm{\kr^{-1/2} \RadFilt(\kr)}{L^{1}(\RR, d\kr)} \norm{\AngFilt(\kt)}{L^{1}(\Sone,d\kt)}
  \, \ConstGeo
\end{multline}
Let us also use the convention that
\begin{equation}
  \label{eq:angularFilterNorm1}
  \norm{\AngFilt(\kt)}{L^{1}(d\kt)} = 1 \, .
\end{equation}

Using the fact that $d(x, \arcsMinusOne{\theta}{\alpha}{j} ) \geq d(x, \arcs{\theta}{\alpha} )$, we find that \eqref{eq:22} implies:
\begin{multline}
  \label{eq:23}
  \contrastMin \sqrt{\frac{\pi}{2 \curvemax}} \inf_{r \in [-\alpha^{2}/2\curvemin, \alpha^{2}/2\curvemin]} \RadFiltCentered(r) /  \ConstGeo \\
  \geq
    \frac{(4M-1)\ConstTangentNormal  }{2M d(x, \arcs{\theta}{\alpha}) }
   + 2 \norm{\kr^{-1/2} \RadFilt(\kr)}{L^{1}(\RR, d\kr)}
\end{multline}
Thus, if $\norm{\kr^{-1/2} \RadFilt(\kr)}{L^{1}(\RR, d\kr)}$ is sufficiently small, then \eqref{eq:22} will be true whenever $d(x, \arcs{\theta}{\alpha})$ is sufficiently large. In particular, if we make each term on the right side of \eqref{eq:23} smaller than half the left side, this equation will be satisfied.

We summarize this result in the following corollary to Theorem \ref{thm:DirectionalFilterYieldsEdges}, which shows that the directional filter is large only when $d(x,\arcs{\theta}{\alpha})$ is sufficiently small.
\begin{corollary}
  \label{thm:ThresholdingSeparatesArcsFromSmoothParts}
  Suppose that:
    \begin{multline}
      \label{eq:asymptoticRemainderSmallEnoughForEdgeIsolation}
      \norm{\kr^{-1/2} \RadFilt(\kr)}{L^{1}(\RR, d\kr)}
  \leq \frac{1}{4}\contrastMin \sqrt{\frac{\pi}{2 \curvemax}} \inf_{r \in [-\alpha^{2}/2\curvemin, \alpha^{2}/2\curvemin]} \RadFiltCentered(r) / \ConstGeo \, .
      \end{multline}
    Define the surfel location error to be:
    \begin{subequations}
      \label{eq:SurfelResolutionThreshold}
      \begin{multline}
        \label{eq:SurfelResolutionDef}
        \SurfelResolution \\
        \equiv \frac{2}{\contrastMin} \sqrt{\frac{2 \curvemax}{\pi}} \left[\inf_{r \in [-\alpha^{2}/2\curvemin, \alpha^{2}/2\curvemin]} \RadFiltCentered(r) \right]^{-1} \frac{(4M-1)\ConstTangentNormal}{2M}
        \, \ConstGeo \, .
      \end{multline}
    \end{subequations}
    Then whenever
    \begin{equation}
      d(x, \arcs{\theta}{\alpha}) \geq \SurfelResolution,
    \end{equation}
    we have that
    \begin{equation}
      \abs{
        \dfilter{\theta}{\alpha} \left[\sum_{j=0}^{M-1} \ImgCoeff{j} \fstep{j}(x) \right](x)
        } \leq \SurfelThreshold.
    \end{equation}
\end{corollary}
Thus, when $x$ is located a distance at least $\SurfelResolution$ from $\arcs{\theta}{\alpha}$, the filtered image is smaller than $\SurfelThreshold$. On the other hand, at the point $\gamma_{j}(t_{j}(\alpha))$, we know that the filtered image is larger than $\SurfelThreshold$. Thus, we obtain the following thresholding algorithm for locating surfels in the wavefront:

\begin{algorithm}[H]
\label{algo:surfelExtraction}
\SetLine
\KwIn{The image in the Fourier domain, i.e. $\Imgk(\vk)$ and a desired minimal sampling rate $\epsilon$.}
\KwOut{Surfels which approximate the wavefront of $\Img(x)$ in the direction $\theta$.}
\klet $f_{1}(x) := [\dfilter{\theta}{\alpha} \Imgk](x)$.\;
\klet $Z := \left\{ x \in [-1,1]^{2} : \abs{f(x)} \geq \SurfelThreshold \right\}$.\;
Cluster the set $Z$. Any two points are part of the same cluster if they are located a distance $\SurfelResolution$ apart. Let $S$ denote the set of clusters. \;
\klet RESULT := []  \quad (the empty set) \;
\ForEach{$s \in S$}{
  Let the midline of $s$ be the set $\{\textrm{midpoint}((x + \theta \RR) \cap s ) : x \in s\}$. \;
  Sample the midline of $s$ with spacing at least $\epsilon$, calling the result $Q$. \;
  \ForEach{$q \in Q$}{
    Add the surfel $(q, \theta)$ to RESULT.\;
    }
}
\Return RESULT\;
\caption{Surfel Extraction}
\end{algorithm}

\begin{figure}
\setlength{\unitlength}{0.240900pt}
\ifx\plotpoint\undefined\newsavebox{\plotpoint}\fi
\sbox{\plotpoint}{\rule[-0.200pt]{0.400pt}{0.400pt}}%
\includegraphics[scale=0.75, bb=85pt 175pt 200pt 525pt]{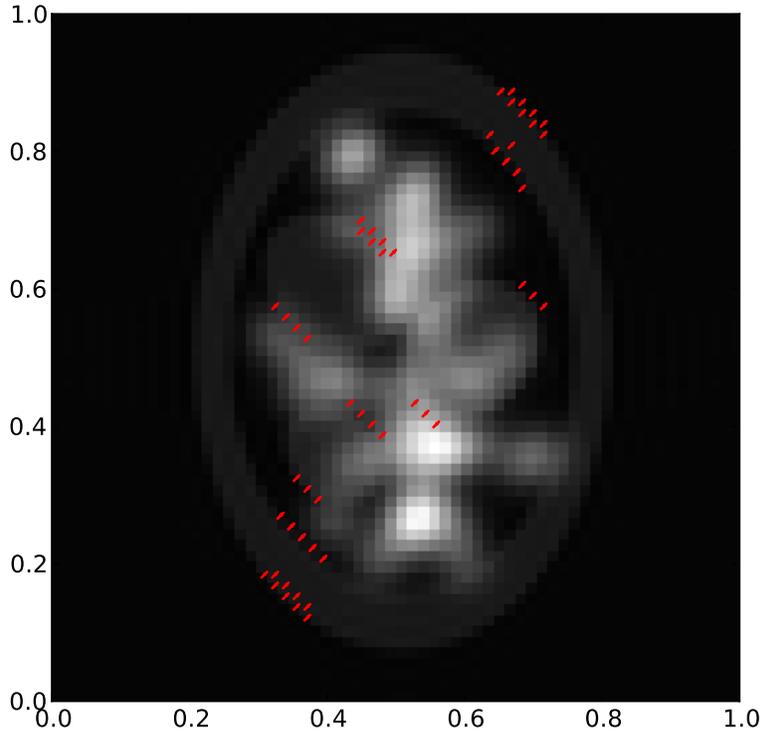}

\caption{An illustration of the result of Algorithm \ref{algo:surfelExtraction}. The red arrows indicate the location and direction of extracted surfels in the direction $\pi/4$ (with $\alpha=\pi/16$). The surfels are overlayed on the image (reconstructed by DFT).}
\label{fig:directionalFilterSurfelExtraction}
\end{figure}

An example of Algorithm \ref{algo:surfelExtraction} applied to the same image as in Section \ref{sec:suboptimalEstimateExample} is shown in Figure \ref{fig:directionalFilterSurfelExtraction}. The algorithm generates no surfel a distance more than $1.0/64$ (i.e., one pixel) away from the actual edge.

We have the following result concerning correctness of Algorithm \ref{algo:surfelExtraction}.

\begin{theorem}
  \label{thm:DirectionalFiltersCorrectlyExtractSurfels}
  Suppose that in addition to \eqref{eq:formOfRadialFilter}, \eqref{eq:geometricConstraintMaxOfFilter}, \eqref{eq:decayOfRadFilterInX}, \eqref{eq:angularFilterNorm1}, \eqref{eq:asymptoticRemainderSmallEnoughForEdgeIsolation}, the following constraint is satisfied:
  \begin{equation}
    \label{eq:SurfelResolutionBiggerThanSeparationDist}
    \SurfelResolution \leq \frac{\delta}{3}
  \end{equation}
  Then for every $(x, \theta)$ in the result of Algorithm \ref{algo:surfelExtraction}, there is a corresponding surfel $(x', \theta')$ in the wavefront of $\Img(x)$ with the property that:
  \begin{subequations}
    \label{eq:SurfelsApproximateWavefront}
    \begin{equation}
      \abs{x - x'} \leq \SurfelResolution
    \end{equation}
    \begin{equation}
      \abs{\theta - \theta'} \leq \alpha
    \end{equation}
  \end{subequations}
  Additionally, for each $j$, at least one surfel output by Algorithm \ref{algo:surfelExtraction} will approximate some surfel $(\gamma_{j}(t_{j}(\lambda)), N_{j}(t_{j}(\lambda)))$ in the arc segment arc segment $\{ \gamma_{j}(t_{j}(\lambda)) : \lambda \in [\theta-\alpha, \theta+\alpha] \}$.
\end{theorem}
\begin{proof}
  Define the arc segment:
  \begin{equation*}
    A_{j} = \{ \gamma_{j}(t_{j}(\lambda)) : \lambda \in [\theta-\alpha, \theta+\alpha] \}
  \end{equation*}

  By Corollary \ref{thm:ThresholdingSeparatesArcsFromSmoothParts}, any point $x \in Z$ is located a distance at most $\SurfelResolution$ away from $\arcs{\theta}{\alpha}$.

Now, consider any segment $s$ of $Z$. For any $x \in s$, there is some point $x'$ in some arc segment $A_{j}$ for which $d(x,x') \leq \SurfelResolution$. We need to show that \emph{all} points $y \in s$ are located at most a distance $\SurfelResolution$ from the \emph{same} arc $A_{j}$. Consider a point $y$ with $d(y, A_{k}) \leq \SurfelResolution$. Suppose also that $y' \in A_{k}$ is a point for which $d(y,y') \leq SurfelResolution$. Then:
\begin{multline*}
  \curvesep \leq d(x', y') \leq d(x', x) + d(x,y) + d(y, y') \leq \SurfelResolution + d(x,y) + \SurfelResolution
\end{multline*}
Subtracting $2 \SurfelResolution$ from both sides and applying \eqref{eq:SurfelResolutionBiggerThanSeparationDist} implies that $d(x,y) \geq \curvesep/3 \geq \SurfelResolution$ and therefore $x$ and $y$ are not in the same segment. Thus, $s$ consists only of points a distance $\SurfelResolution$ from $A_{j}$.

Therefore, any point on the midline of $S$ is located a distance at most $\SurfelResolution$ from $\arcs{\theta}{\alpha}$. By definition, any point on $\arcs{\theta}{\alpha}$ has a normal pointing in some direction in $[\theta-\alpha, \theta+\alpha]$. This proves the Surfels returned by Algorithm \ref{algo:surfelExtraction} accurately approximate surfels in the wavefront of $\Img(x)$.

  To prove that each arc segment has at least one surfel in it, note that by Theorem \ref{thm:DirectionalFilterYieldsEdges}, the point $\gamma_{j}(t_{j}(\theta))$ is contained in $Z$ (c.f. \eqref{eq:DirectionalFilterLargeAtEdge}). This implies that each arc segment $\{ \gamma_{j}(t_{j}(\lambda)) : \lambda \in [\theta-\alpha, \theta+\alpha] \}$ generates at least one segment $s$ of $Z$. There will be at least one sample taken from this segment, which will generate a surfel in the output of Algorithm \ref{algo:surfelExtraction}. Thus, Theorem \ref{thm:DirectionalFiltersCorrectlyExtractSurfels} is proved.
\end{proof}

\section{Segmentation: connecting the surfels is better than connecting the dots}
\label{sec:connectingsurfels}

As we have indicated earlier, one of the reasons for developing a surfel/wavefront extraction
procedure is \emph{segmentation} - by which we mean the reconstruction of the curves
of discontinuity $\gamma_j$ which divide the image into well-defined geometric sub-regions.

The first step in reconstructing the curves is to reconstruct their topology.

\begin{figure}
\setlength{\unitlength}{0.240900pt}
\ifx\plotpoint\undefined\newsavebox{\plotpoint}\fi
\sbox{\plotpoint}{\rule[-0.200pt]{0.400pt}{0.400pt}}%
\includegraphics[scale=0.5]{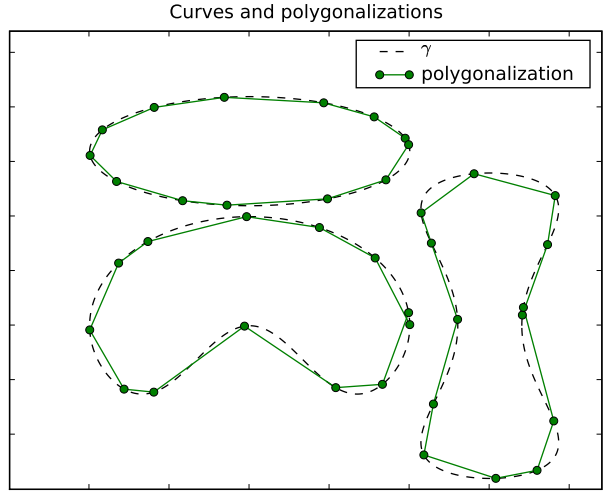}
\caption{A curve and it's polygonalization. }
\label{fig:polygonalization}
\end{figure}

\begin{definition}
  \label{def:polygonalization}
  A polygonalization of a figure $\{\gamma_{j}(t), j = 0,\dots,M-1 \}$ is a planar graph
$\Gamma = (V,E)$ with the property that each vertex $p \in V$ is a point on some $\gamma_{j}(t)$, and each edge connects points which are adjacent samples of some curve $\gamma_{j}(t)$
(see Fig. \ref{fig:polygonalization}).
\end{definition}

\begin{figure}
\setlength{\unitlength}{0.240900pt}
\ifx\plotpoint\undefined\newsavebox{\plotpoint}\fi
\sbox{\plotpoint}{\rule[-0.200pt]{0.400pt}{0.400pt}}%
\centerline{
\includegraphics[scale=0.5]{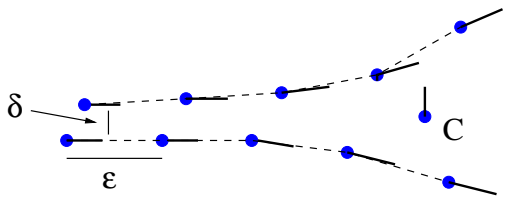}
}
\caption{The advantages of connecting surfels over connecting points are illustrated in this figure. First, the separation between points ($\epsilon$) can be much greater than the separation
between curves ($\delta$). Without information about the tangent, $\epsilon$ and $\delta$
must be of the same order. Using the algorithm of \cite{meleslie:curveppaer}, it is easy to automatically assign points to the correct curve. Note also that points such as $C$ are easy
to filter away when tangent information is available, even in the presence of modest amounts
of noise. }
\label{fig:ctd}
\end{figure}

There is a substantial literature in computational
geometry discussing the task of taking as input a set of unordered points that lie on a set of
curves and returning a polygonalization
\cite{amenta98crust,chengnoise,dey99curve,dey01reconstructing}.
In \cite{meleslie:curveppaer}, we described an algorithm for polygonalization
that uses both point and tangent data (i.e. surfels) and showed that the
method is significantly more robust. It is easier to remove spurious data from the set of
surfels and the sampling requirements are much weaker (see Fig. \ref{fig:ctd}) .

In particular, Theorem 2.8  of  \cite{meleslie:curveppaer} shows that, given a set of points and tangents
(a discrete sampling of the wavefront of $\Img(x)$),
then there is an algorithm that returns the correct polygonalization provided
$\curvesep > 2 \curvemax \epsilon^{2}$ and $\epsilon < (\sqrt{2} \curvemax)^{-1}$, where
$\epsilon$ is the maximal distance between samples on the curve.

If the point and tangent data are corrupted by noise (as they are in practice), then we need to assume a maximal sampling rate as well. Otherwise noise could change the order of
samples on the curve. The following theorem provides technical conditions under which
one can prove that the algorithm is correct when applied to noisy data.

\begin{theorem} \cite[Theorem 3.2]{meleslie:curveppaer}
  \label{thm:noisyPolygonalReconstruction}
  Suppose that Assumptions \ref{ass:boundedCurvature} and \ref{ass:separatedCurves} hold, that noise in the point data is bounded by $\pointNoise$, and that noise in the tangent data is bounded by $\tanNoise$. Suppose further that:
  \begin{subequations}
    \label{eq:noisySeparationConditions}
    \begin{equation}
      \label{eq:noisySeparationDistance}
      \curvesep > 4 \pointNoise + 4 \epsilon \tanNoise + 2.1 \curvemax \epsilon^{2} \, ,
    \end{equation}
    \begin{equation}
      \label{eq:noisyConstraintOnkmaxEpsilon}
      \epsilon < \frac{1}{ \curvemax \sqrt{2}} \, ,
    \end{equation}
  \end{subequations}
  and that adjacent points on a curve are separated by a distance greater than
  $[(1+2^{3/2})(2 \tanNoise \epsilon + \pointNoise)]$.
  Then there is an algorithm that correctly reconstructs the figure.
\end{theorem}

Once the polygonalization of the curve set has been obtained, one can
approximate the geometry with higher order accuracy. This is particularly easy
in the case of surfel data; between each pair of points, cubic Hermite interpolation constructs
a fourth order polynomial in arclength that interpolates the two points and matches the
derivative (tangent) data as well. This achieves fourth order accuracy.

The full segmentation algorithm follows.

\begin{algorithm}[H]
\label{algo:segmentation}
\SetLine
\KwIn{The Fourier transform of the image, $\Imgk(k)$.}
\KwOut{A set of curves approximating the discontinuities of $\Img(x)$.}
\klet S = []\;
\For{$j=1 \ldots A$}{
  \klet $\theta=j \pi/A$\;
  \klet $s =$ result of applying Algorithm \ref{algo:surfelExtraction} (the Surfel Extraction algorithm, see p. \pageref{algo:surfelExtraction}) to $\Imgk(k)$ in the direction $\theta$.\;
  Append $s$ to $S$. \;
  } \tcc{Now $S$ contains surfels pointing in the direction $\theta$ for many values of $\theta$.}\;
  \;
  {\bf Call} the algorithm of \cite{meleslie:curveppaer} to polygonalize $S$. \;
  \Return $C$, where $C$ is the Hermite interpolant of the polygonalization of  $S$.\;
\caption{Segmentation}
\end{algorithm}

This algorithm can be proven ``correct'' in the sense that, for sufficiently large $\kmax$, the algorithm will return a set of curves which are topologically correct.
To do this, we need need to prove that the output of Algorithm \ref{algo:surfelExtraction} meets the requirements of Theorem \ref{thm:noisyPolygonalReconstruction}. This requires verifying that \eqref{eq:noisySeparationConditions} are satisfied. Note that we take Assumptions \ref{ass:boundedCurvature}, \ref{ass:separatedCurves}, \ref{ass:contrastAssumption} and \ref{ass:curvatureBoundedBelow}  as given.

After extracting surfels from the image (as per the loop in lines 2-6 of Algorithm \ref{algo:segmentation}), we find that the error in each surfel's position is bounded by Theorem \ref{thm:DirectionalFiltersCorrectlyExtractSurfels}:
\begin{eqnarray}
  \pointNoise \leq & \SurfelResolution & = O\left(\frac{\sqrt{\kmax^{2} - \ktex^{2}}}{(\kmax-\ktex)^{3/2}} \right)\\
  \tanNoise \leq & \alpha & = O\left( \frac{1}{\sqrt{\kmax + \ktex}}\right)
\end{eqnarray}
Both these quantities are $O(\kmax^{-1/2})$, holding all other factors fixed, and can therefore be made as small as desired.

Note that the angle between adjacent surfels returned by separate applications of Algorithm \ref{algo:surfelExtraction} is at most $2 \alpha + \pi/A$. By Lemma \ref{lemma:arcLengthBetweenAngles}, the separation between two such surfels in arc length is at most $(2\alpha+\pi/A)/\curvemin$. By taking $A=O(\sqrt{\kmax})$ (e.g. $A=\pi/\alpha$), we find that $\epsilon=O(\kmax^{-1/2})$ and thus \eqref{eq:noisyConstraintOnkmaxEpsilon} is satisfied for sufficiently large $\kmax$. To satisfy \eqref{eq:noisySeparationDistance}, observe that:
\begin{equation}
  4 \pointNoise + 4 \epsilon \tanNoise + 2.1 \curvemax \epsilon^{2} \leq O(\kmax^{-1/2} + \kmax^{-1} + \kmax^{-1}) = O(\kmax^{-1/2}) \leq \delta
\end{equation}
For sufficiently large $\kmax$, this implies \eqref{eq:noisySeparationDistance} is satisfied. Thus, we have shown that \eqref{eq:noisyConstraintOnkmaxEpsilon} is satisfied. This implies that for sufficiently large $\kmax$ (holding all other parameters fixed), if Assumptions \ref{ass:boundedCurvature} and \ref{ass:separatedCurves}, \ref{ass:contrastAssumption} and \ref{ass:curvatureBoundedBelow} are satisfied, then Algorithm \ref{algo:segmentation} will successfully segment the image.

We summarize this in a Theorem.
\newpage

\begin{theorem}
  Suppose that Assumption \ref{ass:boundedCurvature}, \ref{ass:separatedCurves}, \ref{ass:contrastAssumption} and \ref{ass:curvatureBoundedBelow} hold. Then for sufficiently large $\ktex$, $\kmax$, Algorithm \ref{algo:segmentation} will successfully approximate the singular support of $\Img(x)$.
\end{theorem}

\begin{figure}
  \includegraphics[scale=0.75, bb = 18pt 180pt 594pt 612pt]{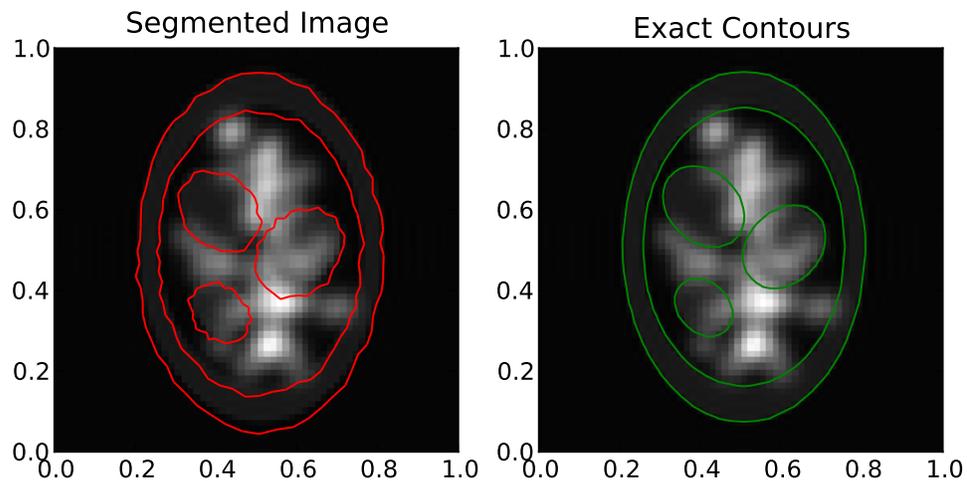}
  \caption{The result of applying the segmentation algorithm to a $64 \times 64$ grid
  of spectral data. Note that the separation of two curves near the center is smaller
  than a pixel, but that wavefront/surfel reconstruction has no difficulty in resolving the
  them.}
  \label{fig:segmentation}
\end{figure}

The result of applying Algorithm \ref{algo:segmentation} to spectral data for our phantom
on a $64 \times 64$ grid is shown in Fig. \ref{fig:segmentation}. The deviation from the exact
result is noticeable, but is on the order of a pixel since we are using low resolution data
to make the nature of the error clear.

\section{Conclusions}

In this paper, we have described a new method for edge detection that can be viewed as
an extension of the method of concentration kernels \cite{tadmor:edgesReview, gelb:gegenbauerReconstructionError,gelbcates:mriSegmentationSpectralData,gelbtadmor:edges,588632}. We use more complex filters in order
to recover information about the wavefront of a two-dimensional image rather than
just its singular set. That is, instead of trying to locate a set of points that lie on curves
of discontinuity, we look for both those point locations and the normal (or tangent) directions there.
This allows us to reconstruct edges more faithfully and robustly, using the algorithm developed previously in \cite{meleslie:curveppaer}. We have focused here on a rigorous mathematical foundation for the method, based on detailed asymptotics and Fourier analysis. Although in this work we require that the curvature of the discontinuities not vanish, this assumption is merely technical. The algorithm works properly even when that assumption is violated, and even for singularities which are not differentiable (see Figure \ref{fig:squareExample} for an example).
A major advantage of extracting surfel information is that one can more easily ``denoise'' the data, as discussed in detail in  \cite{meleslie:curveppaer} and illustrated in Fig. \ref{fig:ctd}.
A number of improvements can still be made, including the
incorporation of nonlinear ``limiters'' to reduce the oscillations produced
in the physical domain from our linear filtering procedure (see, for example,
\cite{588632}).

Recovering local information about a function from partial Fourier data is a rather subtle issue, as demonstrated by Pinsky \cite{pinsky:fourierInversionFails} who showed that spherical partial Fourier integrals do not converge pointwise to the characteristic  function of the unit ball in $\RR^3$. His analysis suggests that radial variations of the concentration method may not converge either (though of course appropriately filtered versions will).

%%In the present context, our analysis shows that we can get good localization of the wavefront in %%the normal direction to the curve, but not in the tangential direction. We expect that this will be %%sufficient for most applications of edge detection.

A limitation of the method described here is that we have assumed the image consists of a globally smooth function superimposed on a set of piecewise constant functions. Extensions of our method to more general piecewise smooth functions will be reported at a later date, as will its application to magnetic resonance imaging.

\appendix

\section{Proof of Proposition \ref{prop:LargeKAsymptoticsOfImage}}
\label{sec:stationaryPhaseCalculations}

To prove Proposition \ref{prop:LargeKAsymptoticsOfImage}, we will require some
results concerning the asymptotics of integrals of the form \eqref{eq:exactTransformOfStepFunction} near a point of stationary phase.

\begin{lemma}
  \label{lemma:stationaryPhaseCurveBoundary}
  Consider a curve $\gamma(t)$, proceeding at unit speed, for $t \in [0,L]$. Suppose $k \cdot \gamma'(0)= 0$, $k \cdot \gamma'(t) \neq 0$ for $t \in (0,L]$, and
  $k \cdot \gamma''(t)$ does not vanish. Let $\ktv = \vk / \abs{\vk}$.

  Then
  \begin{subequations}
    \begin{multline}
      \label{eq:curvedLineSegments}
      \frac{1}{i \abs{k}^{2} } \int_{0}^{L} e^{i k \cdot \gamma(t)}
      \kp \cdot \gamma'(t) dt \\
      =
      \frac{1}{\abs{k}^{3/2}} e^{i k \cdot \gamma(0)} \frac{i^{-1/2} }{2 (\ktv \cdot \gamma''(0)/2)^{1/2}} \gamma(1/2, -i \abs{k} \beta)
      + \frac{R(\vk)}{\abs{k}^{2}}\\
      = \frac{1}{\abs{k}^{3/2}} e^{i k \cdot \gamma(0)} \frac{i^{-1/2} }{2 (\ktv \cdot \gamma''(0)/2)^{1/2}} \Gamma(1/2) \\
      + \frac{R(\vk)}{\abs{k}^{2}} - \abs{k}^{-3/2} e^{i k \cdot \gamma(0)} \frac{i^{-1/2} \Gamma(1/2,-i \abs{k} \beta) }{2 (\ktv \cdot \gamma''(0)/2)^{1/2}}
    \end{multline}
    The remainder $R(\vk)$ is bounded by:
    \begin{equation}
      \abs{R(\vk)} \leq 1 + \frac{3}{4} \frac{ \norm{\ktv \cdot \gamma'''(t)}{L^{\infty}}}{ \abs{\ktv \cdot \gamma''(0)}} \frac{L}{\ktv \cdot \gamma'(L)} + \abs{ \frac{\ktv^{\perp} \cdot \gamma'(L)}{\ktv \cdot \gamma'(L)} - \frac{1 }{\ktv \cdot \gamma''(0) L} }
    \end{equation}
  \end{subequations}
  Here, $\gamma(a,z)$ and $\Gamma(a,z)$ are the incomplete Gamma functions \cite[Chapt. 6.5, p.p.269]{abramowitz:handbookmathfunctions} and $\beta = d \cdot \gamma(L) - d \cdot \gamma(0)$.
\end{lemma}

\begin{proof}
  This is a standard application of stationary phase
  \cite[Section 3.13]{olver:asymptotics}.

  Recalling that , $\ktv = \vk / \abs{k}$, the unit direction of $k$, let us define the variable $v = \ktv \cdot \gamma(t) - \ktv \cdot \gamma(0)$ and $\beta = \ktv \cdot \gamma(L)-\ktv \cdot \gamma(0)$,
  so that $dv = \ktv \cdot \gamma'(t) dt$.
Since $v = \ktv \cdot \gamma(t) - \ktv \cdot \gamma(0) = \ktv \cdot \gamma''(0) t^{2}/2 + O(t^{3})$,
a straightforward calculation shows that
  \begin{equation*}
    \frac{\ktv^{\perp} \cdot \gamma'(t(v))}{\ktv \cdot \gamma'(t(v))} \sim \frac{v^{1/2-1}}{2 (\ktv \cdot \gamma''(0) / 2)^{1/2}} \textrm{  as } v \rightarrow 0.
  \end{equation*}
Thus,
  \begin{multline}
\label{eq:6}
    e^{-i k \cdot \gamma(0)}\int_{0}^{L} e^{i k \cdot \gamma(t)} \kp \cdot \gamma'(t) dt = \abs{k} \int_{0}^{\beta} e^{i \abs{k} v} \frac{\ktv^{\perp} \cdot \gamma'(t(v))}{\ktv \cdot \gamma'(t(v))} dv \\
    =
    \abs{k} \int_{0}^{\beta} e^{i \abs{k} v} \frac{v^{-1/2}}{2 (\ktv \cdot \gamma''(0) / 2)^{1/2}}dv\\
    +
    \abs{k} \int_{0}^{\beta} e^{i \abs{k} v} \left( \frac{\ktv^{\perp} \cdot \gamma'(t(v))}{\ktv \cdot \gamma'(t(v))} - \frac{v^{-1/2}}{2 (\ktv \cdot \gamma''(0) / 2!)^{1/2}} \right) dv \\
    = \abs{k}^{1/2} \frac{i^{1/2} }{2 (\ktv \cdot \gamma''(0)/2)^{1/2}} \gamma(1/2, -i \abs{k} \beta)\\
    +
    \abs{k} \int_{0}^{\beta} e^{i \abs{k} v} \left( \frac{\ktv^{\perp} \cdot \gamma'(t(v))}{\ktv \cdot \gamma'(t(v))} - \frac{v^{-1/2}}{2 (\ktv \cdot \gamma''(0) / 2!)^{1/2}} \right) dv
  \end{multline}
  We must now bound the remainder, the last line of \eqref{eq:6}. This can be done via Theorem 12.3 of \cite[p.p. 99]{olver:asymptotics}, which states that the integral is bounded by the total variation norm of the integrand plus the value at the endpoints, i.e.:
  \begin{multline}
    \label{eq:7}
    \abs{k} \int_{0}^{\beta} e^{i \abs{k} v} \left( \frac{\ktv^{\perp} \cdot \gamma'(t)}{\ktv \cdot \gamma'(t)} - \frac{1 }{\ktv \cdot \gamma''(0) t} \right) dv \\
    = -i e^{i k \cdot \gamma(L)} \left( \frac{\ktv^{\perp} \cdot \gamma'(L)}{\ktv \cdot \gamma'(L)} - \frac{1 }{\ktv \cdot \gamma''(0) L} \right) + ie^{i \vk \cdot \gamma(0)} \left( \frac{\ktv^{\perp} \cdot \gamma''(0)}{\ktv \cdot \gamma''(0)} \right)  + E
  \end{multline}
  with
  \begin{equation*}
    \abs{E} \leq \norm{\left( \frac{\ktv^{\perp} \cdot \gamma'(t)}{\ktv \cdot \gamma'(t)} - \frac{1 }{\ktv \cdot \gamma''(0) t} \right) }{TV} \\
  \end{equation*}
To compute the total variation norm, first use the expansion:
\begin{multline}
  \label{eq:11}
  \left( \frac{\ktv^{\perp} \cdot \gamma'(t)}{d \cdot \gamma'(t)} - \frac{1 }{\ktv \cdot \gamma''(0) t}  \right) = \frac{\sqrt{1-(\ktv \cdot \gamma'(t))^{2}} }{\ktv \cdot \gamma'(t) } - \frac{1 }{\ktv \cdot \gamma''(0) t} \\
  = \frac{1 + (\sqrt{1-(\ktv \cdot \gamma'(t))^{2}} - 1) }{\ktv \cdot \gamma'(t) } - \frac{1 }{\ktv \cdot \gamma''(0) t} \\
  = \frac{ \ktv \cdot \gamma''(0) t - \ktv \cdot \gamma'(t) }{ \ktv \cdot \gamma'(t) \ktv \cdot \gamma''(0) t} + \frac{ \sqrt{1-(\ktv \cdot \gamma'(t))^{2}} - 1 }{ \ktv \cdot \gamma'(t) }
\end{multline}
on $[0,\beta]$.

It is an exercise in elementary calculus to show that $\norm{ f(z) }{TV([0,1])} = 1$ where $f(z)=(\sqrt{1-z^{2}}-1)/z$; this, combined with the fact that $\ktv \cdot \gamma'(t) \in [0,1]$ and $\ktv \cdot \gamma'(t)$ is monotonically increasing, shows that the last term of \eqref{eq:11} has total variation less than $1$.

To bound the first term, we begin by using Taylor's theorem:
  \begin{equation*}
    \ktv \cdot \gamma'(t) - \ktv \cdot \gamma''(0) t \equiv R_{2}(t)  =  \int_{0}^{t} (1/2)\ktv \cdot\gamma'''(z) (z-t) dz
  \end{equation*}
  Define $R(t) = t^{-1} R_{2}(t)$. Then we can write:
  \begin{multline}
    \label{eq:12}
    \norm{ \frac{ \ktv \cdot \gamma''(0) t - \ktv \cdot \gamma'(t) }{ \ktv \cdot \gamma'(t) \ktv \cdot \gamma''(0) t} }{TV} \leq \norm{ \frac{d}{dt} \frac{ \ktv \cdot \gamma''(0) t - \ktv \cdot \gamma'(t) }{ \ktv \cdot \gamma'(t) \ktv \cdot \gamma''(0) t} }{L^{1}} \\
    = \frac{1}{\abs{\ktv \cdot \gamma''(0)}} \int_{0}^{L} \frac{d}{dt} \frac{R(t)}{\ktv \cdot \gamma'(t) } dt = \frac{1}{\abs{\ktv \cdot \gamma''(0)}} \int_{0}^{L} \frac{\ktv \cdot \gamma'(t)  R'(t) - R(t) \ktv \cdot \gamma''(t) }{ (\ktv \cdot \gamma'(t))^{2} } dt
  \end{multline}
  Note that:
  \begin{equation}
    \label{eq:13}
    R'(t) = - \int_{0}^{t} (1/2) \ktv \cdot \gamma'''(z)z/t^{2} dz
  \end{equation}
  Substituting \eqref{eq:13} into \eqref{eq:12} yields:
  \begin{multline}
    \eqref{eq:12} = \frac{1}{2 \abs{\ktv \cdot \gamma''(0)}} \int_{0}^{L} \abs{ \int_{0}^{t} \ktv \cdot \gamma'''(z) \frac{- \ktv \cdot \gamma'(t)  z/t^{2} - (z/t-1) \ktv \cdot \gamma''(t) }{ (\ktv \cdot \gamma'(t))^{2} } } dz dt \\
    = \frac{1}{2 \abs{\ktv \cdot \gamma''(0)}} \int_{0}^{L} \int_{0}^{t} \abs{ \ktv \cdot \gamma'''(z) \frac{- \ktv \cdot \gamma'(t)  z - (zt-t^{2}) \ktv \cdot \gamma''(t) }{ (\ktv \cdot \gamma'(t))^{2} t^{2} } } dz dt \\
    = \frac{1}{2 \abs{\ktv \cdot \gamma''(0)}} \int_{0}^{L} \int_{z}^{L} \abs{ \ktv \cdot \gamma'''(z) \frac{- \ktv \cdot \gamma'(t)  z - (zt-t^{2}) \ktv \cdot \gamma''(t) }{ (\ktv \cdot \gamma'(t))^{2} t^{2} } }dt dz \\
    \leq \frac{ \norm{\ktv \cdot \gamma'''(t)}{L^{\infty}}}{2 \abs{\ktv \cdot \gamma''(0)}}  \int_{0}^{L} \int_{z}^{L} \abs{ \frac{-z}{\ktv \cdot \gamma'(t) t^{2}} - \frac{z \ktv \cdot \gamma''(t) }{ (\ktv \cdot \gamma'(t))^{2} t } + \frac{\ktv \cdot \gamma''(t)}{(\ktv \cdot \gamma'(t))^{2}}}dt dz \\
    \leq \frac{ \norm{\ktv \cdot \gamma'''(t)}{L^{\infty}}}{2 \abs{\ktv \cdot \gamma''(0)}}  \int_{0}^{L} \int_{z}^{L} \abs{ \frac{d}{dt}\left( \frac{z}{\ktv \cdot \gamma'(t) t}\right)  + \frac{\ktv \cdot \gamma''(t)}{(\ktv \cdot \gamma'(t))^{2}}}dt dz \\
    \leq \frac{ \norm{\ktv \cdot \gamma'''(t)}{L^{\infty}}}{2 \abs{\ktv \cdot \gamma''(0)}}
    \int_{0}^{L} \abs{ \frac{-z}{\ktv \cdot \gamma'(L)L} +\frac{z}{\ktv \cdot \gamma'(z)z}  + \frac{1}{\ktv \cdot \gamma'(L)} - \frac{1}{\ktv \cdot \gamma'(z)} } dz \\
    = \frac{ \norm{\ktv \cdot \gamma'''(t)}{L^{\infty}}}{2 \abs{\ktv \cdot \gamma''(0)}} \int_{0}^{L} \frac{z}{\ktv \cdot \gamma'(L)L} + \frac{1}{\ktv \cdot \gamma'(L)} dz \\
    \leq \frac{ \norm{\ktv \cdot \gamma'''(t)}{L^{\infty}}}{2 \abs{\ktv \cdot \gamma''(0)}} \left( \frac{L}{2 \ktv \cdot \gamma'(L)} + \frac{L}{\ktv \cdot \gamma'(L)}\right)\\
    \leq \frac{3}{4} \frac{ \norm{\ktv \cdot \gamma'''(t)}{L^{\infty}}}{ \abs{\ktv \cdot \gamma''(0)}} \frac{L}{\ktv \cdot \gamma'(L)}
  \end{multline}
  Thus, we have the bound that:
  \begin{equation}
     \norm{ \left( \frac{\ktv^{\perp} \cdot \gamma'(t)}{\ktv \cdot \gamma'(t)} - \frac{1 }{\ktv \cdot \gamma''(0) t}  \right) }{TV} \leq 1 + \frac{3}{4} \frac{ \norm{\ktv \cdot \gamma'''(t)}{L^{\infty}}}{ \abs{\ktv \cdot \gamma''(0)}} \frac{L}{\ktv \cdot \gamma'(L)}
  \end{equation}
  This implies that:
  \begin{equation}
    \eqref{eq:7} = -i e^{i k \cdot \gamma(L)} \left( \frac{\ktv^{\perp} \cdot \gamma'(L)}{\ktv \cdot \gamma'(L)} - \frac{1 }{\ktv \cdot \gamma''(0) L} \right) + ie^{i k \cdot \gamma(0)} \left( \frac{\ktv^{\perp} \cdot \gamma''(0)}{\ktv \cdot \gamma''(0)} \right)  + E
  \end{equation}
  with
  \begin{equation*}
    \abs{E} \leq 1 + \frac{3}{4} \frac{ \norm{\ktv \cdot \gamma'''(t)}{L^{\infty}}}{ \abs{\ktv \cdot \gamma''(0)}} \frac{L}{\ktv \cdot \gamma'(L)}
  \end{equation*}
  Since $d$ is the normal vector to $\gamma(t)$ at $t=0$, and $\ktv^{\perp}$ is tangent to it, we find that $\ktv^{\perp} \cdot \gamma''(0) = 0$ (using the fact that $\gamma''(0)=\kappa(0) N(0)$).
  Thus, we find that:
  \begin{multline}
    \label{eq:1}
    \int_{0}^{L} e^{i k \cdot \gamma(t)} \kp \cdot \gamma'(t) dt =
    \abs{k}^{1/2} \frac{i^{1/2} }{2 (\ktv \cdot \gamma''(0)/2)^{1/2}} \gamma(1/2, -i \abs{k} \beta) \\
    + i e^{i k \cdot \gamma(0)} \left( \frac{\ktv^{\perp} \cdot \gamma''(0)}{\ktv \cdot \gamma''(0)} \right) + E' = \abs{k}^{1/2} \frac{i^{1/2} }{2 (\ktv \cdot \gamma''(0)/2)^{1/2}} \gamma(1/2, -i \abs{k} \beta) + E'
  \end{multline}
  with
  \begin{equation}
    \abs{E'} \leq 1 + \frac{3}{4} \frac{ \norm{\ktv \cdot \gamma'''(t)}{L^{\infty}}}{ \abs{\ktv \cdot \gamma''(0)}} \frac{L}{\ktv \cdot \gamma'(L)} + \abs{ \frac{\ktv^{\perp} \cdot \gamma'(L)}{\ktv \cdot \gamma'(L)} - \frac{1 }{\ktv \cdot \gamma''(0) L} }
  \end{equation}
  Multiplying \eqref{eq:1} by $1/(i \abs{k}^{2})$ yields the result we seek.
\end{proof}

We also need the following geometric result:

\begin{lemma}
  \label{lemma:arcLengthBetweenAngles}
  Consider two normal vectors $k_{1}$ and $k_{2}$ on a curve $\gamma_{j}(t)$
  with an angle $\theta$ between them. Then:
  \begin{equation}
    \label{eq:arcLengthBetweenAngles}
    \theta / \curvemin \geq \abs{t_{j}(k_{1}) - t_{j}(k_{2})} \geq \abs{\theta} / \curvemax
  \end{equation}
\end{lemma}
\begin{proof}
  The angle changes most quickly if $\gamma_{j}(t)$ is a circle with minimal radius of curvature. The arc length along such a curve is $\abs{\theta} R \geq \abs{\theta} / \curvemax$. The arc length changes most quickly (w.r.t angle) along a circle with maximal radius of curvature, i.e. a circle of radius $1/\curvemin$.
\end{proof}

The proof of Proposition \ref{prop:LargeKAsymptoticsOfImage} basically requires us to applying Lemma \ref{lemma:stationaryPhaseCurveBoundary} to \eqref{eq:exactTransformOfStepFunction}.

\begin{proofof}{Proposition  \ref{prop:LargeKAsymptoticsOfImage}}
  Note that $k \cdot \gamma_{j}'(t_{j}(\vk)) = 0$. Thus, we may apply Lemma \ref{lemma:stationaryPhaseCurveBoundary} to the curve $\gamma(t) = \gamma_{j}(t_{j}(\vk)+t)$ along the segment $t \in [0,t_{k}(k^{\perp}) - t_{j}(\vk)]$ and similarly to $\gamma(t)=\gamma_{j}(t_{j}(\vk)-t)$. This will give us an expansion over the sections of the curve where $k \cdot N_{j}(t) \geq 0$. Repeating the analysis centered at $t_{j}(-k)$ yields an expansion over sections of the curve where $k \cdot N_{j}(t) \leq 0$.

Applying Lemma \ref{lemma:stationaryPhaseCurveBoundary} directly yields the following:
\begin{multline}
  \label{eq:2}
  \frac{1}{i \abs{k}^{2} } \int_{k \cdot N_{j}(t) > 0} e^{i k \cdot \gamma(t)}
  \kp \cdot \gamma'(t) dt  =
  \frac{1}{\abs{k}^{3/2}} e^{i k \cdot \gamma(0)} \frac{i^{-1/2} }{2 (\ktv \cdot \gamma''(t_{j}(\vk))/2)^{1/2}} \Gamma(1/2) \\
      + \frac{R_{1}(\vk)}{\abs{k}^{2}} - \abs{k}^{-3/2} e^{i k \cdot \gamma(0)} \frac{i^{-1/2} \Gamma(1/2,-i \abs{k} \beta_{1}) }{2 (\ktv \cdot \gamma''(t_{j}(\vk))/2)^{1/2}} \\
      + \frac{1}{\abs{k}^{3/2}} e^{i k \cdot \gamma(0)} \frac{i^{-1/2} }{2 (\ktv \cdot -\gamma''(t_{j}(\vk))/2)^{1/2}} \Gamma(1/2) \\
      + \frac{R_{2}(\vk)}{\abs{k}^{2}} - \abs{k}^{-3/2} e^{i k \cdot \gamma(0)} \frac{i^{-1/2} \Gamma(1/2,-i \abs{k} \beta_{2}) }{2 (\ktv \cdot -\gamma''(t_{j}(\vk))/2)^{1/2}} \\
      = e^{i k \cdot \gamma(t_{j}(\vk))} \frac{\sqrt{\pi} }{\abs{k}^{3/2} \sqrt{\kappa_{j}(t_{j}(\vk))} } + \frac{R_{1}(\vk)+R_{2}(\vk)}{\abs{k}^{2}}\\
      + \abs{k}^{-3/2} e^{i k \cdot \gamma(0)} \frac{i^{-1/2} \Gamma(1/2,-i \abs{k} \beta_{1}) }{2 (\ktv \cdot \gamma''(t_{j}(\vk))/2)^{1/2}} + \abs{k}^{-3/2} e^{i k \cdot \gamma(0)} \frac{i^{-1/2} \Gamma(1/2,-i \abs{k} \beta_{2}) }{2 (\ktv \cdot \gamma''(t_{j}(\vk))/2)^{1/2}}
\end{multline}

Here, the remainders $R_{1,2}(\vk)$ are bounded by:
\begin{subequations}
  \label{eq:3}
  \begin{multline}
    \abs{R_{1}(\vk)} \leq 1 + \frac{3}{4} \frac{ \norm{\ktv \cdot \gamma_{j}'''(t)}{L^{\infty}}}{ \abs{\ktv \cdot \gamma''(t_{j}(\vk))}} \frac{(t_{j}(k^{\perp}) - t_{j}(\vk))}{\ktv \cdot \gamma'(t_{j}(k^{\perp}))}\\
    + \abs{ \frac{\ktv^{\perp} \cdot \gamma'(t_{j}(k^{\perp}))}{\ktv \cdot
        \gamma'(t_{j}(k^{\perp}))} - \frac{1 }{\ktv \cdot
        \gamma''(t_{j}(\vk))(t_{j}(k^{\perp}) - t_{j}(\vk))} }
  \end{multline}
  \begin{multline}
    \abs{R_{1}(\vk)} \leq 1 + \frac{3}{4} \frac{ \norm{\ktv \cdot \gamma_{j}'''(t)}{L^{\infty}}}{ \abs{\ktv \cdot \gamma''(t_{j}(\vk))}} \frac{(t_{j}(\vk) - t_{j}(-k^{\perp}))}{\ktv \cdot \gamma'(t_{j}(-k^{\perp}))}\\
    + \abs{ \frac{\ktv^{\perp} \cdot \gamma'(t_{j}(-k^{\perp}))}{\ktv \cdot
        \gamma'(t_{j}(-k^{\perp}))} - \frac{1 }{\ktv \cdot
        \gamma''(t_{j}(\vk))(t_{j}(\vk) - t_{j}(-k^{\perp}))}}
  \end{multline}
\end{subequations}
Since $k \cdot \gamma_{j}'(t_{j}(\vk))=0$ and $\gamma_{j}(t)$ proceeds with unit speed, we find that $\ktv \cdot \gamma_{j}''(t_{j}(\vk)) = \kappa(t_{j}(\vk))$. Note also that $\ktv^{\perp} \cdot \gamma_{j}''(t_{j}(\pm k^{\perp})) = 0$ since $\gamma_{j}''(t_{j}(\pm k^{\perp})) \parallel \pm k^{\perp}$ (and similarly $\ktv \cdot \gamma'(t_{j}(-k^{\perp})) = 1$. Substituting this into \eqref{eq:3}, as well as bounding the curvature below by $\curvemin$ and adding them up yields:
\begin{multline}
  \label{eq:4}
  \abs{R_{1}(\vk)} + \abs{R_{2}(\vk)} \leq 2 + \frac{3}{2} \frac{\norm{\gamma_{j}'''(t)}{L^{\infty}} }{ \curvemin} \frac{ t_{j}(k^{\perp}) - t_{j}(-k^{\perp})}{ 1}\\
  + \abs{ \frac{1}{\kappa_{j}(t_{j}(\vk))(t_{j}(k^{\perp}) - t_{j}(\vk))}} + \abs{ \frac{1}{\kappa_{j}(t_{j}(\vk))(t_{j}(\vk) - t_{j}(-k^{\perp}) )}}
\end{multline}
Applying Lemma \ref{lemma:arcLengthBetweenAngles} shows that $(t_{j}(k^{\perp}) - t_{j}(\vk)) \geq \pi/2 \curvemax$ and similarly $(t_{j}(\vk) - t_{j}(-k^{\perp}) ) \geq \pi/2 \curvemax$. Substituting this into \eqref{eq:4} and observing that $\kappa_{j}(t_{j}(\vk)) > \curvemin$ yields:
\begin{multline}
  \label{eq:8}
  \abs{R_{1}(\vk)} + \abs{R_{2}(\vk)} \leq 2 + \frac{3}{2} \frac{\norm{\gamma_{j}'''(t)}{L^{\infty}} }{ \curvemin} (t_{j}(k^{\perp}) - t_{j}(-k^{\perp})) + \frac{4 \curvemax}{ \pi \curvemin}
\end{multline}
We must also bound
\begin{equation*}
  \abs{
    \abs{k}^{-3/2} e^{i k \cdot \gamma(0)} \frac{i^{-1/2} \Gamma(1/2,-i \abs{k} \beta_{1}) }{2 (\ktv \cdot \gamma''(t_{j}(\vk))/2)^{1/2}} + \abs{k}^{-3/2} e^{i k \cdot \gamma(0)} \frac{i^{-1/2} \Gamma(1/2,-i \abs{k} \beta_{2}) }{2 (\ktv \cdot \gamma''(t_{j}(\vk))/2)^{1/2}}
    }
\end{equation*}
To do this, note that $\beta_{1}=\ktv \cdot \gamma_{j}(t_{j}(k^{\perp})) - \ktv \cdot \gamma_{j}(t_{j}(\vk)) \leq - 1/\curvemax$ (and similarly $\beta_{2} \leq -1/\curvemax$. Thus, $-\abs{k} \beta_{1,2} \geq \abs{\kr}\curvemax$. This can be seen easily by considering circle tangent to $\gamma_{j}(t_{j}(\vk))$ of radius $1/\curvemax$. Using Eq. (6.5.32) from \cite[p.p. 263]{abramowitz:handbookmathfunctions}, combined with the estimate on the remainder stated immediately after Eq. (6.5.32), we observe that $\abs{\Gamma(1/2,-i \abs{k} \beta_{i,j})} \leq \abs{k/ \curvemax}^{-1/2}$. Thus, we find that:
\begin{multline}
  \abs{
    \abs{k}^{-3/2} e^{i k \cdot \gamma(0)} \frac{i^{-1/2} \Gamma(1/2,-i \abs{k} \beta_{1}) }{2 (\ktv \cdot \gamma''(t_{j}(\vk))/2)^{1/2}} + \abs{k}^{-3/2} e^{i k \cdot \gamma(0)} \frac{i^{-1/2} \Gamma(1/2,-i \abs{k} \beta_{2}) }{2 (\ktv \cdot \gamma''(t_{j}(\vk))/2)^{1/2}}
  } \\
  \leq \abs{ \abs{k}^{-3/2}  \frac{(k \curvemax)^{-1/2}}{ 2(\curvemin/2)}} + \abs{\abs{k}^{-3/2}  \frac{(k \curvemax)^{-1/2}}{ 2(\curvemin/2)}}
  \leq \frac{1}{\abs{k}^{2}} \sqrt{ \frac{2 \curvemax}{ \curvemin} }
\end{multline}
Repeating this analysis for the part of the curve where $k \cdot N_{j}(t) < 0$ yields the following:
\begin{equation}
  \fstep{j}(\vk) = \frac{ e^{i k \cdot \gamma(t_{j}(\vk))}}{\abs{k}^{3/2}} \frac{\sqrt{\pi} }{\sqrt{\kappa_{j}(t_{j}(\vk))} } + \frac{ e^{i k \cdot \gamma(t_{j}(-k))}}{\abs{k}^{3/2}} \frac{\sqrt{\pi} }{ \sqrt{\kappa_{j}(t_{j}(-k))} } + \frac{E_{j}(\vk)}{\abs{k}^{2}}
\end{equation}
where
\begin{equation}
  \abs{E_{j}(\vk)} \leq 4 + 3\frac{\norm{\gamma_{j}'''(t)}{L^{\infty}} }{ \curvemin} \textrm{arclength}(\gamma_{j}) + \frac{8 \curvemax}{ \pi \curvemin} + 2 \sqrt{ \frac{2 \curvemax}{ \curvemin} }
\end{equation}

Adding this result up over $j=0\ldots M-1$ and bounding $\ImgCoeff{j}$ by $\contrastMax$ yields the desired result.
\end{proofof}

\section{Proof of Theorem \ref{thm:DirectionalFilterYieldsEdges}: Leading Order Asymptotics}
\label{sec:DirectionalFiltersLeadingOrderAsymptotocs}

We first show that the directional filters behave properly when applied to the leading order asymptotic terms of $\Imgk(\vk)$.

For this, we need to prove two facts: a) that directional filters, after being applied to the image, decay away from the points $\gamma_{j}(t_{j}(\theta))$, with $\theta$ the direction of the filter
and b) that the filters yield spikes at or near the points $\gamma_{j}(t_{j}(\theta))$.

The basis for our calculation is the following Lemma, which allows us to write the directional filter applied to the leading order term of \eqref{eq:asymptoticExpansionOfStepFunction} as an integral over the curve $\gamma_{j}(t)$.

We consider only a directional filter oriented in the direction $\theta=0$. Results for
other directions can be obtained by rotation.

\begin{lemma}
  \label{lemma:angularInverseFourierOfFilter}
  Let $\AngFilt(\kt)$ be supported on the interval $[-\alpha,\alpha]$, let $\RadFilt(\kr) \in L^{1}(\RR, d\kr)$ and let $\RadFilt(\kr)/\kr \in L^{1}(\RR, d\kr)$. Then:
  \begin{multline}
    \label{eq:angularInverseFourierOfFilter}
    \int_{0}^{\infty} \int_{-\alpha}^{\alpha} e^{-i k \cdot x}\frac{ e^{i k \cdot \gamma_{j}(t_{j}(\vk))}}{\abs{k}^{3/2}} \frac{\sqrt{\pi} }{\sqrt{\kappa_{j}(t_{j}(\vk))} } \AngFilt(\kt) \kr^{1/2} \RadFilt(\kr) d\kt \kr d\kr \\
    = \int_{-\alpha}^{\alpha} \int_{0}^{\infty} e^{i \kr N_{j}(t_{j}(\kt)) \cdot (\gamma_{j}(t_{j}(\kt))-x)} \frac{\sqrt{\pi} }{\sqrt{\kappa_{j}(t_{j}(\kt))} } \AngFilt(\kt) \RadFilt(\kr) d\kr d\kt \\
    = \int_{-\alpha}^{\alpha} \frac{\sqrt{\pi} }{\sqrt{\kappa_{j}(t_{j}(\kt))} } \AngFilt(\kt) \check{\RadFilt}(N_{j}(t_{j}(\kt)) \cdot [\gamma_{j}(t_{j}(\kt))-x]) d\kt \\
    = \sqrt{\pi} \int_{t_{j}(-\alpha)}^{t_{j}(\alpha)} \AngFilt(\kt(t)) \sqrt{\kappa_{j}(t)} \check{\RadFilt}(N_{j}(t) \cdot [\gamma_{j}(t)-x]) dt
  \end{multline}
\end{lemma}

\begin{proof}
  \begin{multline}
    \label{eq:14}
    \int_{0}^{\infty} \int_{-\alpha}^{\alpha} e^{-i k \cdot x}\frac{ e^{i k \cdot \gamma_{j}(t_{j}(\vk))}}{\abs{k}^{3/2}} \frac{\sqrt{\pi} }{\sqrt{\kappa_{j}(t_{j}(\vk))} } \AngFilt(\kt) \kr^{1/2} \RadFilt(\kr) d\kt \kr d\kr \\
    = \int_{-\alpha}^{\alpha} \int_{0}^{\infty} e^{i \kr N_{j}(t_{j}(\kt)) \cdot (\gamma_{j}(t_{j}(\kt))-x)} \frac{\sqrt{\pi} }{\sqrt{\kappa_{j}(t_{j}(\kt))} } \AngFilt(\kt) \RadFilt(\kr) d\kr d\kt
  \end{multline}
  Note that the inner integral of the last line of \eqref{eq:14} is merely the inverse Fourier transform of $\RadFilt(\kr)$ evaluated at the point $r=N_{j}(t_{j}(\kt)) \cdot [\gamma_{j}(t_{j}(\kt))-x]$. Thus,
  \begin{multline*}
    \eqref{eq:14}
    = \int_{-\alpha}^{\alpha} \frac{\sqrt{\pi} }{\sqrt{\kappa_{j}(t_{j}(\kt))} } \AngFilt(\kt) \check{\RadFilt}(N_{j}(t_{j}(\kt)) \cdot [\gamma_{j}(t_{j}(\kt))-x]) d\kt \\
    = \sqrt{\pi} \int_{t_{j}(-\alpha)}^{t_{j}(\alpha)} \AngFilt(\kt(t)) \sqrt{\kappa_{j}(t)} \check{\RadFilt}(N_{j}(t) \cdot [\gamma_{j}(t)-x]) dt \, ,
  \end{multline*}
  completing the proof.
\end{proof}

\begin{proposition}
  \label{prop:boundOnFilteredTangentNormal}
  Let $\AngFilt(\kt)$ be supported on the interval $[-\alpha,\alpha]$, let $\RadFilt(\kr) \in L^{1}(\RR, d\kr)$ and also $\RadFilt(\kr)/\kr \in L^{1}(\RR, d\kr)$.

  Define the smallest normal and tangent distances ($D_{N}$ and $D_{T}$ respectively) as:
  \begin{subequations}
    \begin{eqnarray}
      D_{N} & = & \inf_{t \in [t_{j}(-\alpha),t_{j}(\alpha)]}N_{j}(t) \cdot [\gamma_{j}(t) - x]\\
      D_{T} & = & \inf_{t \in [t_{j}(-\alpha),t_{j}(\alpha)]}\gamma_{j}'(t) \cdot [\gamma_{j}(t) - x].
    \end{eqnarray}
  \end{subequations}
  \begin{subequations}
    Then we have the following bound on the action of the filter in
    the normal directions:
    \begin{multline}
      \label{eq:boundOnFilteredNormal}
      \abs{ \int_{0}^{\infty} \int_{-\alpha}^{\alpha} e^{-i k \cdot x}\frac{ e^{i k \cdot \gamma_{j}(t_{j}(\vk))}}{\abs{k}^{3/2}} \frac{\sqrt{\pi} }{\sqrt{\kappa_{j}(t_{j}(\vk))} } \AngFilt(\kt) \kr^{1/2} \RadFilt(\kr) d\kt \kr d\kr } \\
      \leq \sqrt{\frac{\pi}{\curvemin}}
      \norm{\AngFilt(\kt)}{L^{1}(\Sone,d\kt)} \RadFiltXSup(D_{N})
    \end{multline}
    We also have a weaker bound in the tangential direction:
    \begin{multline}
      \label{eq:boundOnFilteredTangent}
      \abs{ \int_{0}^{\infty} \int_{-\alpha}^{\alpha} e^{-i k \cdot x}\frac{ e^{i k \cdot \gamma_{j}(t_{j}(\vk))}}{\abs{k}^{3/2}} \frac{\sqrt{\pi} }{\sqrt{\kappa_{j}(t_{j}(\vk))} } \AngFilt(\kt) \kr^{1/2} \RadFilt(\kr) d\kt \kr d\kr } \\
      \leq \frac{2 \sqrt{\pi} \norm{\RadFilt(\kr)/\kr}{L^{1}(\RR, d\kr)}}{D_{T}} \left( \frac{\norm{\AngFilt'(\kt)}{L^{1}}}{\curvemin^{1/2}} + \frac{\norm{\AngFilt(\kt)}{L^{1}}\norm{\gamma_{j}'''(t)}{L^{\infty}}}{2 \curvemin^{5/2}}\right)\\
      + \frac{\sqrt{\pi} \norm{\AngFilt(\kt)}{L^{1}}
        \norm{\check{\RadFilt}(z)(z+1)}{L^{\infty}}}{D_{T}^{2}}
    \end{multline}
  \end{subequations}
  Finally,
  \begin{equation}
    \label{eq:infOfNormalDist}
    D_{N} \geq \cos(\alpha) \abs{N_{j}(t_{j}(0))
          \cdot (\gamma(t_{j}(0)) - x)} - \frac{\alpha^{2}}{\curvemin^{2}} \left[
          \curvemin + \frac{\curvemax}{2} +
          \frac{\alpha}{\curvemin^{3}} \frac{\norm{\gamma'''(t)}{L^{\infty}}}{6} \right]
  \end{equation}
\end{proposition}

\begin{proof}

  The result \eqref{eq:boundOnFilteredNormal} follows from the second to last line of \eqref{eq:angularInverseFourierOfFilter} and the fact that
  \begin{equation*}
    \abs{\check{\RadFilt}(N_{j}(t_{j}(\kt)) \cdot [\gamma_{j}(t_{j}(\kt))-x])} \leq \RadFiltXSup(D_{N}).
  \end{equation*}

To prove \eqref{eq:infOfNormalDist}, we must bound $N_{j}(t_{j}(\kt)) \cdot [\gamma(t_{j}(\kt)) - x]$. Let $t_{0}=t_{j}(0)$, and consider the Taylor expansion (to second order) of $\gamma_{j}(t)$.
  \begin{equation*}
    N_{j}(t) \cdot \left[\gamma_{j}(t_{0}) + \gamma_{j}'(t_{0})(t-t_{0}) + \frac{\gamma_{j}''(t_{0})}{2}(t-t_{0})^{2} + \textrm{remainder} \right]
  \end{equation*}
By Taylor's theorem, the remainder is bounded by $\abs{k \cdot R(t)} \leq \norm{\gamma'''(t)}{L^{\infty}} \frac{\alpha^{3}}{6 \kappa_{j}(t)}$. The first order term is:
  \begin{equation*}
    \abs{N_{j}(\kt) \cdot \gamma'(t_{0}) } = \abs { \cos(\kt + \pi/2) } \leq  \abs{\sin(\kt - 0) }
    \leq \abs{\kt} \leq \alpha
  \end{equation*}
  We use also the fact that $\abs{t-t_{0}} \leq \alpha/\curvemin$.
  Thus, we obtain:
  \begin{multline}
    \label{eq:differenceInNormalDistance}
    \abs{N_{j}(t) \cdot (\gamma(t)-x) - N_{j}(t) \cdot( \gamma(t_{0}) - x)} \leq \frac{\alpha^{2}}{\curvemin} + \frac{\curvemax}{2} \frac{\alpha^{2}}{\curvemin^{2}}  + \frac{\norm{\gamma'''(t)}{L^{\infty}}}{6}\abs{t-t_{0}}^{3} \\
    \leq \alpha^{2}\left(\frac{1}{\curvemin} + \frac{\curvemax}{2 \curvemin^{2}}\right) + \frac{\norm{\gamma'''(t)}{L^{\infty}}}{6} \frac{\alpha^{3}}{\curvemin^{3}}
  \end{multline}
  Note also that $N_{j}(t) \cdot ( \gamma(t_{0}) - x) \geq \cos(\alpha) N_{j}(t_{0}) \cdot ( \gamma(t_{0}) - x)$. We therefore find that:
  \begin{multline}
    \label{eq:17}
    \abs{N_{j}(t) \cdot [\gamma(t) - x] } \geq\\
    \cos(\alpha) \abs{N_{j}(t_{0}) \cdot (\gamma(t_{0}) - x)} - \left[ \alpha^{2}\left(\frac{1}{\curvemin} + \frac{\curvemax}{2 \curvemin^{2}}\right) + \frac{\norm{\gamma'''(t)}{L^{\infty}}}{6} \frac{\alpha^{3}}{\curvemin^{3}} \right]
  \end{multline}
  Taking the inf over the angles in $[-\alpha,\alpha]$ yields the result we seek.

  For \eqref{eq:boundOnFilteredTangent}, note that $\partial_{t} N_{j}(t) \cdot [\gamma_{j}(t) - x] = \kappa_{j}(t) \gamma_{j}'(t) \cdot [\gamma_{j}(t) - x]$. We can then multiply and divide the integrand of \eqref{eq:angularInverseFourierOfFilter} by this and then integrate by parts to obtain:
  \begin{multline}
    \label{eq:15}
    \eqref{eq:angularInverseFourierOfFilter}\\
    =\sqrt{\pi} \int_{t_{j}(-\alpha)}^{t_{j}(\alpha)} \frac{
      \AngFilt(\kt(t)) \sqrt{\kappa_{j}(t)} }{
      \kappa_{j}(t) \gamma_{j}'(t) \cdot [\gamma_{j}(t) - x]
      }
      \check{\RadFilt}(N_{j}(t) \cdot [\gamma(t)-x])  \kappa_{j}(t) \gamma_{j}'(t) \cdot [\gamma_{j}(t) - x] dt \\
      = - \sqrt{\pi} \int_{t_{j}(-\alpha)}^{t_{j}(\alpha)} \check{\RadFilt}(N_{j}(t) \cdot [\gamma(t)-x])
      \partial_{t} \frac{
        \AngFilt(\kt(t))
      }{
        \sqrt{\kappa_{j}(t)} \gamma_{j}'(t) \cdot [\gamma_{j}(t) - x]
      } dt \\
      = - \sqrt{\pi} \int_{t_{j}(-\alpha)}^{t_{j}(\alpha)} \check{\RadFilt}(N_{j}(t) \cdot [\gamma(t)-x])
      \times \Bigg(
      \frac{ \AngFilt'(\kt(t)) \kt'(t) } { \sqrt{\kappa_{j}(t)} \gamma_{j}'(t) \cdot [\gamma_{j}(t) - x]} \\
      + \frac{ \AngFilt(\kt(t))  \kappa_{j}'(t) }{ 2(\kappa_{j}(t))^{3/2} \gamma_{j}'(t) \cdot [\gamma_{j}(t) - x] }
      + \frac{
        \AngFilt(\kt(t)) \left( \gamma''(t) \cdot [\gamma_{j}(t) - x] + \gamma'(t) \cdot \gamma'(t)\right)
      }{
        \kappa_{j}(t) (\gamma_{j}'(t) \cdot [\gamma_{j}(t) - x] )^{2}
      }
      \Bigg)dt \\
      = - \sqrt{\pi} \int_{t_{j}(-\alpha)}^{t_{j}(\alpha)} \check{\RadFilt}(N_{j}(t) \cdot [\gamma(t)-x])
      \times \Bigg(
      \frac{ \AngFilt'(\kt(t)) \kappa_{j}(t) } { \sqrt{\kappa_{j}(t)} \gamma_{j}'(t) \cdot [\gamma_{j}(t) - x]} \\
      + \frac{ \AngFilt(\kt(t))  \kappa_{j}'(t) }{ 2(\kappa_{j}(t))^{3/2} \gamma_{j}'(t) \cdot [\gamma_{j}(t) - x] }
      + \frac{
        \AngFilt(\kt(t)) \left( \kappa_{j}(t) N_{j}(t) \cdot [\gamma_{j}(t) - x] + 1 \right)
      }{
        \kappa_{j}(t) (\gamma_{j}'(t) \cdot [\gamma_{j}(t) - x] )^{2}
      }
      \Bigg)dt
  \end{multline}
  We can further simplify this to:
  \begin{multline}
    \label{eq:16}
    \eqref{eq:15} = - \sqrt{\pi} \int_{t_{j}(-\alpha)}^{t_{j}(\alpha)} \check{\RadFilt}(N_{j}(t) \cdot [\gamma(t)-x])
      \times \Bigg(
      \frac{ \AngFilt'(\kt(t)) \sqrt{\kappa_{j}(t)} } { \gamma_{j}'(t) \cdot [\gamma_{j}(t) - x]} \\
      + \frac{ \AngFilt(\kt(t))  \kappa_{j}'(t) }{ 2(\kappa_{j}(t))^{3/2} \gamma_{j}'(t) \cdot [\gamma_{j}(t) - x] }
      + \frac{
        \AngFilt(\kt(t)) \left( N_{j}(t) \cdot [\gamma_{j}(t) - x] + 1 \right)
      }{
        (\gamma_{j}'(t) \cdot [\gamma_{j}(t) - x] )^{2}
      }
      \Bigg)dt \\
      = \sqrt{\pi} \int_{-\alpha}^{\alpha} \check{\RadFilt}(N_{j}(t_{j}(\kt)) \cdot [\gamma(t_{j}(\kt))-x])
      \times \Bigg(
      \frac{ \AngFilt'(\kt(t_{j}(\kt))) \sqrt{\kappa_{j}(t_{j}(\kt))} } { \gamma_{j}'(t_{j}(\kt)) \cdot [\gamma_{j}(t_{j}(\kt)) - x]} \\
      + \frac{ \AngFilt(\kt(t_{j}(\kt)))  \kappa_{j}'(t_{j}(\kt)) }{ 2(\kappa_{j}(t_{j}(\kt)))^{3/2} \gamma_{j}'(t_{j}(\kt)) \cdot [\gamma_{j}(t_{j}(\kt)) - x] } \\
      +
      \frac{
        \AngFilt(\kt(t_{j}(\kt))) \left( N_{j}(t_{j}(\kt)) \cdot [\gamma_{j}(t_{j}(\kt)) - x] + 1 \right)
      }{
        (\gamma_{j}'(t_{j}(\kt)) \cdot [\gamma_{j}(t_{j}(\kt)) - x] )^{2}
      }
      \Bigg) \frac{d\kt}{\kappa_{j}(t_{j}(\kt))}
  \end{multline}
  Note that $\check{F}(z) = \int^{z} \check{f}(z') dz' = \int e^{-i \kr z} \RadFilt(\kr)(i\kr)^{-1} d\kr$, and thus $\abs{\check{F}(z)} \leq \norm{\RadFilt(\kr)/\kr}{L^{1}(\RR^{1})}$. Noting also that $\gamma_{j}'''(t) \cdot N_{j}(t) = \kappa_{j}'(t)$ (differentiate the formula $\gamma_{j}''(t)=\kappa_{j}(t) N_{j}(t)$, use the Frenet-Serret formula and dot product with $N_{j}(t)$), we find that $\abs{\kappa_{j}'(t)} \leq \abs{\gamma_{j}'''(t)}$. Combining these facts and using the definition of $D_{T}$, we obtain:
  \begin{multline}
    \abs{ \eqref{eq:16} } \\ \leq
    \sqrt{\pi} \int_{-\alpha}^{\alpha} \norm{\RadFilt(\kr)/\kr}{L^{1}(\RR, d\kr)}
    \Bigg(
      \frac{ \AngFilt'(\kt(t_{j}(\kt)))  } {\sqrt{\kappa_{j}(t_{j}(\kt))} D_{T}}
      + \frac{ \AngFilt(\kt(t_{j}(\kt)))  \norm{\gamma_{j}'''(t)}{L^{\infty}} }{ 2(\kappa_{j}(t_{j}(\kt)))^{5/2} D_{T} } \Bigg) \\
      + \check{F}(N_{j}(t_{j}(\kt)) \cdot [\gamma(t_{j}(\kt))-x]) (N_{j}(t_{j}(\kt)) \cdot [\gamma(t_{j}(\kt))-x] + 1)
      \frac{
        \AngFilt(\kt(t_{j}(\kt)))
      }{
        D_{T}^{2}
      }
      \Bigg) d\kt \\
      \leq     \sqrt{\pi} \int_{-\alpha}^{\alpha} \norm{\RadFilt(\kr)/\kr}{L^{1}(\RR, d\kr)}
    \Bigg(
      \frac{ \AngFilt'(\kt(t_{j}(\kt)))  } {\sqrt{\kappa_{j}(t_{j}(\kt))} D_{T}}
      + \frac{ \AngFilt(\kt(t_{j}(\kt)))  \norm{\gamma_{j}'''(t)}{L^{\infty}} }{ 2(\kappa_{j}(t_{j}(\kt)))^{5/2} D_{T} } \Bigg) \\
      + \norm{\check{F}(z)(z+1)}{L^{\infty}}
      \frac{
        \AngFilt(\kt(t_{j}(\kt)))
      }{
        D_{T}^{2}
      }
      d\kt \\
      \leq \frac{2 \sqrt{\pi} \norm{\RadFilt(\kr)/\kr}{L^{1}(\RR, d\kr)}}{D_{T}} \left( \frac{\norm{\AngFilt'(\kt)}{L^{1}(\Sone,d\kt)}}{\curvemin^{1/2}} + \frac{\norm{\AngFilt(\kt)}{L^{1}}\norm{\gamma_{j}'''(t)}{L^{\infty}}}{2 \curvemin^{5/2}}\right)\\
      + \frac{\sqrt{\pi} \norm{\AngFilt(\kt)}{L^{1}} \norm{\check{F}(z)(z+1)}{L^{\infty}}}{D_{T}^{2}}
  \end{multline}
  This yields \eqref{eq:boundOnFilteredTangent}, and completes the proof of Proposition \ref{prop:boundOnFilteredTangentNormal}.
\end{proof}

\begin{lemma}
  \label{lemma:normalDistanceVsAngle}
  Let $\gamma(t)$ be a curve moving at unit speed and having non-vanishing curvature, with unit tangent $T(t)$ and normal $N(t)$. Then:
  \begin{equation}
    \label{eq:normalDistanceVsAngle}
    \abs{N(t) \cdot (\gamma(t)-\gamma(t_{0}))} \leq \frac{\theta^{2}(t)}{2 \curvemin}
  \end{equation}
  Here, $\theta(t) = T(t) \angle T(t_{0})$.
\end{lemma}

\begin{proof}
  Letting $\theta(t)$ be the angle of the tangent, we find that $\theta'(t)=\kappa(t)$. Note that non-vanishing curvature implies that $t(\theta)$ and $\theta(t)$ are both functions (at least for small $t$ and $\theta$). Note first that:
  \begin{equation*}
    \frac{d}{d\theta} (\gamma(t)-\gamma(t_{0})) = \gamma'(t) \frac{dt}{d\theta} = \frac{T(t)}{\kappa(t)}
  \end{equation*}
  Integrating this with respect to $\theta$ shows that $\abs{\gamma(t)-\gamma(t_{0})} \leq \theta/\curvemin$. Now compute:
  \begin{multline*}
    \frac{d}{d\theta} N(t) \cdot (\gamma(t)-\gamma(t_{0})) = \frac{1}{\kappa(t)} \frac{d}{dt} N(t) \cdot (\gamma(t)-\gamma(t_{0})) \\
    = \frac{1}{\kappa(t)} (-\kappa(t) T(t)) \cdot (\gamma(t) - \gamma(t_{0})) + \frac{1}{\kappa(t)} N(t) \cdot T(t) = T(t) \cdot (\gamma(t) - \gamma(t_{0}))
  \end{multline*}
  This implies that:
  \begin{equation*}
    \abs{\frac{d}{d\theta} N(t) \cdot (\gamma(t)-\gamma(t_{0}))} \leq \abs{T(t) \cdot (\gamma(t) - \gamma(t_{0}))} \leq \frac{\theta}{\curvemin}
  \end{equation*}
  Integrating with respect to $\theta$ yields the result we seek.
\end{proof}

\begin{proposition}
  \label{prop:maxOfFilter}
  Suppose that $\AngFilt(\kt)$ is smooth and compactly supported on $[-\alpha,\alpha]$. Let $\RadFiltCenteredX(r)$ be a function symmetric about $r=0$, and strictly positive on the interval $[-\alpha^{2}/2\curvemin, \alpha^{2}/2 \curvemin]$. Suppose also that $\RadFilt(\kr)$ satisfies \eqref{eq:formOfRadialFilter}. Assume also that \eqref{eq:geometricConstraintMaxOfFilter} is satisfied.

    Then at the point $\gamma_{j}(t_{j}(0))$, we have that:
    \begin{multline}
      \int_{-\infty}^{\infty} \int_{-\alpha}^{\alpha} e^{-i k \cdot x}\frac{ e^{i k \cdot \gamma_{j}(t_{j}(\vk))}}{\abs{k}^{3/2}} \frac{\sqrt{\pi} }{\sqrt{\kappa_{j}(t_{j}(\vk))} } \AngFilt(\kt) \kr^{1/2} \RadFilt(\kr) d\kt \kr d\kr \\
      \geq \sqrt{\frac{\pi}{2 \curvemax}} \norm{\AngFilt(\kt)}{L^{1}(d\kt)}
      \inf_{r \in [-\alpha^{2}/2\curvemin, \alpha^{2}/2\curvemin]}
      \RadFiltCentered(r)
    \end{multline}

\end{proposition}

\begin{proof}
  Note that:
  \begin{equation}
    \label{eq:10}
    k \cdot (\gamma_{j}(t_{j}(\vk)) - x) = \kr N_{j}(t_{j}(\kt)) \cdot [\gamma_{j}(t_{j}(\kt)) - x]
  \end{equation}
  Following the calculations of \eqref{eq:angularInverseFourierOfFilter} and using \eqref{eq:inverseFourierTransformOfRadialFilter} as well as \eqref{eq:10} we find:
  \begin{multline}
    \label{eq:9}
    \int_{-\infty}^{\infty} \int_{-\alpha}^{\alpha} e^{-i k \cdot x}\frac{ e^{i k \cdot \gamma_{j}(t_{j}(\vk))}}{\abs{k}^{3/2}} \frac{\sqrt{\pi} }{\sqrt{\kappa_{j}(t_{j}(\vk))} } \AngFilt(\kt) \kr^{1/2} \RadFilt(\kr) d\kt \kr d\kr \\
    = \int_{-\alpha}^{\alpha} \frac{\sqrt{\pi} }{\sqrt{\kappa_{j}(t_{j}(\kt))} } \AngFilt(\kt) \check{\RadFilt}(N_{j}(t_{j}(\kt)) \cdot [\gamma_{j}(t_{j}(\kt))-x]) d\kt \\
    = \int_{-\alpha}^{\alpha} \frac{\sqrt{\pi} }{\sqrt{\kappa_{j}(t_{j}(\kt))} } \AngFilt(\kt)  g\Big(N_{j}(t_{j}(\kt)) \cdot [\gamma_{j}(t_{j}(\kt))-x]\Big) \\
    \times \cos \Big( \frac{\kmax + \ktex}{2}N_{j}(t_{j}(\kt)) \cdot [\gamma_{j}(t_{j}(\kt))-x]\Big) d\kt
  \end{multline}
  For $x = \gamma_{j}(t_{j}(0))$ and $\alpha \in [-\alpha,\alpha]$, Lemma \ref{lemma:normalDistanceVsAngle} implies that:
  \begin{equation*}
    \abs{N_{j}(t) \cdot (\gamma(t)-x)} \leq \frac{\alpha^{2}}{2\curvemin} \leq \frac{\pi}{4}
  \end{equation*}
  where the last line follows from Lemma \ref{lemma:normalDistanceVsAngle}. This implies that $\cos({N_{j}(t) \cdot (\gamma(t)-x)}) \geq 1/\sqrt{2}$ and therefore:
  \begin{multline*}
    \eqref{eq:9} \geq \int_{-\alpha}^{\alpha} \frac{\sqrt{\pi} }{\sqrt{\kappa_{j}(t_{j}(\kt))} } \AngFilt(\kt) m \frac{1}{\sqrt{2}} d\kt\\
    \geq \sqrt{\pi/2} \norm{\AngFilt(\kt)}{L^{1}(d\kt)} / \sqrt{\curvemax} \inf_{r \in [-\alpha^{2}/2\curvemin, \alpha^{2}/2\curvemin]} \RadFiltCenteredX(r)
  \end{multline*}
  This is what we wanted to prove.
\end{proof}

\subsection{Proof of Theorem \ref{thm:DirectionalFilterYieldsEdges}: Putting it together}

\label{sec:PuttingDirectionalFilterTheoremTogether}

We now consider the behavior of the filter applied to the entire asymptotic expansion of the image.

\begin{lemma}
  \label{lemma:decayOfDirectionalFilterWhenRadialWindowIsInverseR}
  Suppose that $\RadFilt(\kr)$ satisfies \eqref{eq:decayOfRadFilterInX}. In that case
    \begin{equation}
      \abs {\dfilter{\theta}{\alpha}
          \frac{ e^{i k \cdot \gamma(t_{j}(\vk))}}{\abs{k}^{3/2}}
            \frac{\sqrt{\pi} }{\sqrt{\kappa_{j}(t_{j}(\vk))} } }
          \leq \ConstTangentNormal \inf_{t \in [t_{j}(-\alpha),t_{j}(\alpha)]}
          \frac{1}{\abs{\gamma_{j}(t)-x}},
    \end{equation}
    where $\ConstTangentNormal$ is given by \eqref{eq:defOfConstTangentNormal}.
\end{lemma}

\begin{proof}
  As before, we study the case when $\theta=0$, since the rest can be treated by rotation.

  Recalling \eqref{eq:infOfNormalDist}, we define the quantities:
  \begin{eqnarray*}
    D_{N}^{+} & = & \sup_{s \in [t_{j}(-\alpha),t_{j}(\alpha)]} \inf_{t \in [t_{j}(-\alpha),t_{j}(\alpha)]} \abs{N_{j}(s) \cdot [\gamma_{j}(t) - x]}\\
    D_{T}^{+} & = & \sup_{s \in [t_{j}(-\alpha),t_{j}(\alpha)]}\inf_{t \in [t_{j}(-\alpha),t_{j}(\alpha)]} \abs{T_{j}(s) \cdot [\gamma_{j}(t) - x]}.
    \end{eqnarray*}
    Obviously $D_{N}^{+} \geq D_{N}$ and $D_{T}^{+} \geq D_{T}$. Let $N=N_{j}(s)$ be the normal at which $D_{N}^{+}$ is achieved, and $T=T_{j}(s)$ be the tangent at which $D_{T}^{+}$ is achieved. Let $\vp = \gamma_{j}(t) - x$. Since the angle between $N$ and $T$ is at most $2\alpha$, we find that $\abs{x}^{2} \leq \sec(2\alpha)^{2} \abs{x}^{2}_{(N,T)}$, where $\abs{~\cdot~}_{(N,T)}$ is the $L^{2}$ norm taken in the coordinate system $(N,T)$. Thus, we find:
    \begin{equation*}
      \abs{\vp}^{2} \leq \sec(2\alpha)^{2} \abs{\vp}^{2}_{(N,T)} \leq  2 \sec(2\alpha)^{2} \max \{ \abs{N \cdot \vp}, \abs{T \cdot \vp} \}^{2}
    \end{equation*}
    Thus:
    \begin{multline*}
      \sup_{t \in [t_{j}(-\alpha),t_{j}(\alpha)]} \abs{\gamma_{j}(t) - x} \\
      \leq \sqrt{2}\sec(2\alpha)
      \sup_{t \in [t_{j}(-\alpha),t_{j}(\alpha)]} \max \{ \abs{N \cdot (\gamma_{j}(t) - x)}, \abs{T \cdot (\gamma_{j}(t) - x)} \}^{2} \\
      = \sqrt{2} \sec(2\alpha) \max \{ \sup_{t \in [t_{j}(-\alpha),t_{j}(\alpha)]} \abs{N \cdot (\gamma_{j}(t) - x)}, \sup_{t \in [t_{j}(-\alpha),t_{j}(\alpha)]} \abs {T\cdot (\gamma_{j}(t) - x)} \} \\
      = \sqrt{2} \sec(2\alpha) \max\{ D_{N}^{+}, D_{T}^{+} \}
    \end{multline*}
    This implies that:
    \begin{equation}
      \label{eq:20}
      \inf \frac{1}{\abs{\gamma_{j}(t)-x}} \geq \frac{\cos(2\alpha)}{\sqrt{2}} \min\left\{ \frac{1}{D_{N}^{+}}, \frac{1}{D_{T}^{+}} \right\} \geq \frac{\cos(2\alpha)}{\sqrt{2}} \min\left\{ \frac{1}{D_{N}}, \frac{1}{D_{T}} \right\}
    \end{equation}

    Now, substituting \eqref{eq:decayOfRadFilterInX} into \eqref{eq:boundOnFilteredNormal} we find:
  \begin{equation}
    \label{eq:18}
    \abs{\dfilter{\alpha}{\theta}
      \frac{ e^{i k \cdot \gamma_{j}(t_{j}(\vk))}}{\abs{k}^{3/2}} \frac{\sqrt{\pi} }{\sqrt{\kappa_{j}(t_{j}(\vk))} } }
    \leq
    \frac{\ConstRadFiltDecay \sqrt{\pi / \curvemin} \norm{\AngFilt(\kt)}{L^{1}(\Sone,d\kt)}}{D_{N}}
  \end{equation}
  Using \eqref{eq:boundOnFilteredTangent}, we find:
  \begin{multline}
    \label{eq:21}
    \abs{ \dfilter{\alpha}{\theta} \frac{ e^{i k \cdot \gamma_{j}(t_{j}(\vk))}}{\abs{k}^{3/2}} \frac{\sqrt{\pi} }{\sqrt{\kappa_{j}(t_{j}(\vk))} }  } \\
      \leq \frac{2 \sqrt{\pi} \norm{\RadFilt(\kr)/\kr}{L^{1}(\RR, d\kr)}}{ D_{T}} \left( \frac{\norm{\AngFilt'(\kt)}{L^{1}}}{\curvemin^{1/2}} + \frac{\norm{\AngFilt(\kt)}{L^{1}}\norm{\gamma_{j}'''(t)}{L^{\infty}}}{2 \curvemin^{5/2}}\right)\\
      + \frac{\sqrt{\pi} \norm{\AngFilt(\kt)}{L^{1}} \norm{\check{\RadFilt}(z)(z+1)}{L^{\infty}}}{D_{T}^{2}} \\
      \leq \frac{\sqrt{\pi}}{D_{T} \sqrt{\curvemin}} \Bigg(
      2 \norm{\RadFilt(\kr)/\kr}{L^{1}(\RR, d\kr)} \left[\norm{\AngFilt'(\kt)}{L^{1}}
        + \norm{\AngFilt(\kt)}{L^{1}}  \frac{\norm{\gamma_{j}'''(t)}{L^{\infty}}}{2 \curvemin^{2}} \right]\\
        + \sqrt{\curvemin} \norm{\check{\RadFilt}(z)(z+1)}{L^{\infty}}
        \Bigg)
  \end{multline}
  Taking the min of \eqref{eq:18} and \eqref{eq:21}, and using \eqref{eq:20} we find:
  \begin{multline}
    \abs{ \dfilter{\alpha}{\theta} \frac{ e^{i k \cdot \gamma_{j}(t_{j}(\vk))}}{\abs{k}^{3/2}} \frac{\sqrt{\pi} }{\sqrt{\kappa_{j}(t_{j}(\vk))} }  } \leq \frac{\sqrt{\pi}}{\sqrt{\curvemin} }  \min\left\{ \frac{1}{D_{N}}, \frac{1}{D_{T}} \right\} \\
    \times \max\Bigg\{
      \ConstRadFiltDecay \norm{\AngFilt(\kt)}{L^{1}(\Sone,d\kt)},\\
\Bigg(
      2 \norm{\RadFilt(\kr)/\kr}{L^{1}(\RR, d\kr)} \left[\norm{\AngFilt'(\kt)}{L^{1}}
        + \norm{\AngFilt(\kt)}{L^{1}}  \frac{\norm{\gamma_{j}'''(t)}{L^{\infty}}}{2 \curvemin^{2}} \right]\\
        + \sqrt{\curvemin} \norm{\check{\RadFilt}(z)(z+1)}{L^{\infty}}
        \Bigg)
      \Bigg\} \\
      \leq \left( \inf \frac{1}{\abs{\gamma_{j}(t)-x}} \right) \frac{\sqrt{2\pi}}{\sqrt{\curvemin} \cos(2\alpha)} \\
      \max\Bigg\{
      \ConstRadFiltDecay \norm{\AngFilt(\kt)}{L^{1}(\Sone,d\kt)},\\
\Bigg(
      2 \norm{\RadFilt(\kr)/\kr}{L^{1}(\RR, d\kr)} \left[\norm{\AngFilt'(\kt)}{L^{1}}
        + \norm{\AngFilt(\kt)}{L^{1}}  \frac{\norm{\gamma_{j}'''(t)}{L^{\infty}}}{2 \curvemin^{2}} \right]\\
        + \sqrt{\curvemin} \norm{\check{\RadFilt}(z)(z+1)}{L^{\infty}}
        \Bigg)
      \Bigg\}
  \end{multline}
  This is what we wanted to show.
\end{proof}

Recalling Definition \ref{def:arcsOfSpecificAngle}, we have the following result:
\begin{proposition}
  \label{prop:filterBoundedAwayFromArcs}
  Let $x$ be a point not on any of the curves $\gamma_{i}(t)$. Then (with $\ConstTangentNormal$ as in Lemma \ref{lemma:decayOfDirectionalFilterWhenRadialWindowIsInverseR}):
  \begin{multline}
    \abs {\dfilter{\theta}{\alpha}
      \left[\sum_{j=0}^{M-1} \ImgCoeff{j} \fstep{j}(x) \right] }
      \leq
      \frac{\ConstTangentNormal}{d(x, \arcs{\theta}{\alpha} ) }\\
      + \norm{\kr^{-1/2} \RadFilt(\kr)}{L^{1}(\RR, d\kr)} \norm{\AngFilt(\kt)}{L^{1}(\Sone,d\kt)} \,
      \ConstGeo \, ,
  \end{multline}
  where $\ConstGeo$ is defined in (\ref{eq:asymptoticExpansionOfStepFunctionRemainder2}).
\end{proposition}

\begin{proof}
  Begin by using \eqref{eq:asymptoticExpansionOfStepFunction}:
  \begin{multline*}
    \abs{\dfilter{\theta}{\alpha} \sum_{j=0}^{M-1} \ImgCoeff{j} \fstep{j}(\vk)}\\
    \leq \abs{\dfilter{\theta}{\alpha} \sum_{j=0}^{M-1} \ImgCoeff{j} \left[
      \frac{ e^{i k \cdot \gamma(t_{j}(\vk))}}{\abs{k}^{3/2}} \frac{\sqrt{\pi} }{\sqrt{\kappa_{j}(t_{j}(\vk))} } + \frac{ e^{i k \cdot \gamma(t_{j}(-k))}}{\abs{k}^{3/2}} \frac{\sqrt{\pi} }{ \sqrt{\kappa_{j}(t_{j}(-k))} }
      \right]}
      + \abs{\dfilter{\theta}{\alpha}\frac{E(\vk)}{\abs{k}^{2}}} \\
      \leq \sum_{j=0}^{M-1} \ImgCoeff{j} \left[
        \abs{\dfilter{\theta}{\alpha} \frac{ e^{i k \cdot \gamma(t_{j}(\vk))}}{\abs{k}^{3/2}} \frac{\sqrt{\pi} }{\sqrt{\kappa_{j}(t_{j}(\vk))} }} + \abs{\dfilter{\theta}{\alpha} \frac{ e^{i k \cdot \gamma(t_{j}(-k))}}{\abs{k}^{3/2}} \frac{\sqrt{\pi} }{ \sqrt{\kappa_{j}(t_{j}(-k))} }}
      \right] \\
      + \abs{\dfilter{\theta}{\alpha}\frac{E(\vk)}{\abs{k}^{2}}}
      \leq \ConstTangentNormal \\
      \times \sum_{j=0}^{M-1} \ImgCoeff{j} \left[ \inf_{t \in [t_{j}(-\alpha+\theta),t_{j}(\alpha + \theta)]}
      \frac{1}{\abs{\gamma_{j}(t)-x}} + \inf_{t \in [t_{j}(-\alpha+\theta + \pi),t_{j}(\alpha + \theta + \pi)]}
      \frac{1}{\abs{\gamma_{j}(t)-x}} \right] \\
    + \abs{\dfilter{\theta}{\alpha}\frac{E(\vk)}{\abs{k}^{2}}} \\
    \leq \frac{\ConstTangentNormal 2M }{d(x, \arcs{\theta}{\alpha} ) } \max_{j} \ImgCoeff{j}  + \abs{\dfilter{\theta}{\alpha}\frac{E(\vk)}{\abs{k}^{2}}} \\
    \leq \frac{\ConstTangentNormal 2M }{d(x, \arcs{\theta}{\alpha} ) } \contrastMax + \abs{\dfilter{\theta}{\alpha}\frac{E(\vk)}{\abs{k}^{2}}}
  \end{multline*}
  To control the last term, recall the bound $\ConstGeo$ in\eqref{eq:asymptoticExpansionOfStepFunctionRemainder2} on $E(\vk)$:
  \begin{multline*}
    \abs{\dfilter{\theta}{\alpha}\frac{E(\vk)}{\abs{k}^{2}}} = \abs{\int e^{-i k \cdot x} \AngFilt(\kt) \kr^{1/2} \RadFilt(\kr) \frac{E(\vk)}{\kr^{2}} d\kt \kr d\kr} \\
    \leq  \int \int \abs{\frac{\RadFilt(\kr)}{\kr^{1/2}}  \AngFilt(\kt) E(\vk) } d\kt d \kr \\
    \leq \int \abs{\frac{\RadFilt(\kr)}{\kr^{1/2}}} d\kr \int \ConstGeo \abs{\AngFilt(\kt) } d\kt \\
    \leq  \norm{\kr^{-1/2} \RadFilt(\kr)}{L^{1}(\RR, d\kr)} \norm{\AngFilt(\kt)}{L^{1}(\Sone,d\kt)} \,
    \ConstGeo.
  \end{multline*}
\end{proof}

\begin{proposition}
  \label{prop:filterLargeAtRightPlace}
  We have the following estimate at the point $x=\gamma_{j}(t_{j}(\theta))$.
  \begin{multline}
    [\theta \cdot N_{j}(t)] \dfilter{\theta}{\alpha} \left[
      \sum_{j=0}^{M-1} \ImgCoeff{j} \fstep{j}(x)
    \right] \geq
    \contrastMin \sqrt{\frac{\pi}{2 \curvemax}} \norm{\AngFilt(\kt)}{L^{1}(d\kt)}
    \inf_{r \in [-\alpha^{2}/2\curvemin, \alpha^{2}/2\curvemin]}
    \RadFiltCentered(r)\\
    - \left[ \frac{\ConstTangentNormal (2M-1) }{2M \, d(x,
        \arcsMinusOne{\theta}{\alpha}{j} ) }
      + \norm{\kr^{-1/2} \RadFilt(\kr)}{L^{1}(\RR, d\kr)} \norm{\AngFilt(\kt)}{L^{1}(\Sone,d\kt)}\right] \,
      \ConstGeo \, .
    \end{multline}
  \end{proposition}
  \begin{proof}
    By Proposition \ref{prop:maxOfFilter}, we have that:
    \begin{multline*}
      \dfilter{\theta}{\alpha} \contrastMin \frac{ e^{i k \cdot \gamma_{j}(t_{j}(\vk))}}{\abs{k}^{3/2}} \frac{\sqrt{\pi} }{\sqrt{\kappa_{j}(t_{j}(\vk))} } 1_{\kt \in [-\alpha,\alpha]}(\kt) \\
      \geq \sqrt{\frac{\pi}{2 \curvemax}} \contrastMin \norm{\AngFilt(\kt)}{L^{1}(d\kt)}
      \inf_{r \in [-\alpha^{2}/2\curvemin, \alpha^{2}/2\curvemin]}
      \RadFiltCentered(r)
    \end{multline*}
    Applying the same reasoning as in the proof of Proposition \ref{prop:filterBoundedAwayFromArcs}, we find that:
    \begin{multline*}
      \abs{ \dfilter{\theta}{\alpha} \left[\sum_{i \neq j} \ImgCoeff{j} \fstep{i}(x) \right] }
      \leq \\
      \left[ \frac{\ConstTangentNormal (2M-2) }{d(x,\arcsMinusOne{\theta}{\alpha}{j} ) }
        + \norm{\kr^{-1/2} \RadFilt(\kr)}{L^{1}(\RR, d\kr)} \norm{\AngFilt(\kt)}{L^{1}(\Sone,d\kt)}\right]
        \, \ConstGeo \, .
         \end{multline*}
    The second triangle inequality yields the result we seek.
  \end{proof}

  \begin{proofof}{Theorem \ref{thm:DirectionalFilterYieldsEdges}}
    Proposition \ref{prop:filterBoundedAwayFromArcs} proves \eqref{eq:DirectionalFilterYieldsEdges}.

    Note that Assumption \ref{ass:separatedCurves} implies that $d(\gamma_{j}(t_{j}(\theta)), \arcsMinusOne{\theta}{\alpha}{j}) \geq \curvesep$. Substituting this into Proposition \ref{prop:filterLargeAtRightPlace} proves \eqref{eq:DirectionalFilterLargeAtEdge}, and thus completes the proof of Theorem \ref{thm:DirectionalFilterYieldsEdges}.
  \end{proofof}

\nocite{gelb:gegenbauerReconstructionError}
\nocite{gelb:gegenbauerReconstructionNonuniform}
\nocite{1066589}
\nocite{gelb:gibbsRingingTissueBoundary}

\bibliography{../stucchio}
\bibliographystyle{plain}

\end{document}